\documentclass[12pt,twoside]{article}
\usepackage{amssymb,amsmath,amsthm,latexsym}
\usepackage{amscd}
\usepackage{graphicx}

\let\dis\displaystyle
\let\p\partial

\def\Rho{P}
\def\nc{{>\!\!\!\!<}}

\def\D{\mathcal{D}}
\def\B{\mathcal{B}}
\def\A{\mathcal{A}}
\def\Rep{{\rm Rep\,}}
\def\Hom{{\rm Hom}}
\def\Ob{{\rm Ob\,}}

\def\ind{{\rm ind\,}}

\def\Mor{{\rm Mor\,}}

\newtheorem{lemma}{Lemma}

\newtheorem{proposition}{Proposition}

\newtheorem{remark}{Remark}
\newtheorem{example}{Example}

\newcommand{\pbf}[1]{\mbox{\mathversion{bold}$#1$}}

\begin{document}

\begin{center}
\LARGE \bf The norm of a relation, separating \\ functions and
representations\\ of marked quivers
\end{center}

\begin{center}
\Large\bf L.A. Nazarova, A.V. Roiter
\end{center}

\noindent Institute of Mathematics of NAS of Ukraine, 3
Tereshchenkivska Str., 01601 Kyiv 4, Ukraine\\ E-mail:
roiter@imath.kiev.ua

\bigskip

{\small \noindent Numerical functions, which characterize Dynkin
schemes, Coxeter graphs and tame marked quivers, are considered.}

\bigskip

It is well known that Dynkin schemes (and extended ones)
are widespread in modern mathematics.
In essence, it is also known that at least a part of these schemes can
be characterized by some numerical functions.
However, sometimes these functions arise as well,
although Dynkin schemes cannot be introduced in a natural way.

The authors found such a situation studying the representations of
dyadic and tryadic sets \cite{4,5}, where the functions
$\mu(n_1,n_2,n_3)=n_1n_2+n_1n_3+n_2n_3+n_1n_2n_3$ and the
conditions $\mu<4$ and $\mu=4$ of finite representativity and
tameness repeatedly appear.

On the other hand,
in \cite{1,2} it was introduced the norm $\|R\|\in[0,1]$ of a binary relation $R$
having a substantial interest at least for partially ordered sets
in connection with their representations~\cite{3}.

If $(S,\leq)$ is a poset, then letting $\Rho(S)=\|\leq\|^{-1}$
we obtain the criteria of finite representativity and tameness of
posets in the form $\Rho(S)<4$ and $\Rho(S)=4$ respectively.
If $S_t$ is a disjoint union of $t$ chains of orders $n_1,\dots,n_t$,
whose elements are pairwise incomparable,
then
$\dis\Rho(S_t)=\sum\limits_{i=1}^{t}\left(1+{n_i-1\over
n_i+1}\right)$ (see lemmas 3,~4).

The trivial transformation for $t=3$ shows that $\Rho(S_t)\leq
4$ if and only if $\mu(n_1-1,n_2-1,n_3-1)\leq 4$.

Identifying $S_t$ with the vector $(n_1,\dots,n_t)$
it is possible to consider a numerical function $\rho$
(of any number of variables~$t$)
defined on ${\mathbb N}$ by the formula above
and extended in a natural way if some variables equal~$\infty$.
For $t=3$ and finite values of variables
this function in fact coincides (up to a factor and an addend, both constant)
with the function $n_1^{-1}+n_2^{-1}+n_3^{-1}.$
The latter plays an important role in the theory of Coxeter groups, systems of roots, etc.
and goes back to M\"obius~\cite{6}
$(\displaystyle \rho (n_1,n_2,n_3)= 6-2\sum\limits_{i=1}^3 (n_i+1)^{-1})$.

In sections~1 and 2 we consider the function $\Rho$ defined for
arbitrary finite set $S$ with the binary relation $R$. It follows
from the definition that if $S'\subseteq S,$ then
$\Rho(S',R)\leq\Rho(S,R)$. Studying the $\Rho$-faithful posets,
i.e. such posets $S$ that $\Rho(S',R)<\Rho(S,R)$ for $S'\subset
S,$ essentially clarifies the structure of the posets $K_1-K_5$
and $N_0-N_5$ in the criteria of finite representativity and
tameness of posets \cite{7,8}.

In sections~3 and 4 the function $\rho$ is naturally considered
as a function of several variables defined on $({\mathbb N}\sqcup\infty),$
and Dynkin schemes, extended Dynkin schemes and Coxeter graphs
are characterized in these terms.

In sections~6--8 the function $\rho$ is applied for studying of
representations of marked quivers introduced by the authors
\cite{4,9}. Marked quivers are a generalization of Gabriel quivers
\cite{10} and contain most of the known matrix problems~\cite{11}.
Here it turns out that the results formulated in the terms of Dynkin schemes
can be often reformulated in terms  of the function~$\rho,$
but in several cases the functions $\rho$ and $\mu$ are necessary,
whereas using of Dynkin schemes is impossible or difficult.

\section{The norm of a binary relation}

Let $R$ be an arbitrary binary relation on a finite set
$S=\{1,\dots,n\}$. $R$ can be considered as a function
of two variables:
define $R(i,j)=1$ for $(i,j)\in R$ and $R(i,j)=0$ otherwise.
Let $f_R$ be a real quadratic form of $n$ variables defined via
$f_R(x_1,\dots,x_n)=\sum\limits_{i,j=1}^{n}R(i,j)x_ix_j,$ and
let $\|R\|$ be the least value taken by $f_R$ on the
standard simplex, i.e. on the set of vectors $(x_1,\dots,x_n),$
where $x_i'$s are nonnegative real numbers and $\sum\limits_{i=1}^{n}x_i=1$.
We call the number $\|R\|$ the {\it norm} of the relation $R$.
This value is reached by an elementary calculus theorem and
can be calculated by usual rules that implies, in particular,
that $\|R\|$ is a rational number.

We call a vector $\bar{x}=(x_1,\dots,x_n)$ (in general, not unique)
such that $f_R(\bar x)=\|R\|$ the {\it minimal vector} for $(S,R)$.

It is easy to see, that $0\leq\|R\|\leq 1$, moreover, $\|R\|=0$ iff
$R$ is not reflexive and $\|R\|=1$ if (and only if) $R$ is a
complete relation (i.e.  $iRj$ for all $i,j\in S$). Further we
will assume that all considering relations are reflexive
($iRi$ for all $i\in S$)\footnote{The replacement of minimum by maximum
in the definition of the norm does not seem meaningful;
if we define such a norm $\overline{\|\cdot\|},$ it is easy to see that
$\overline{\|R\|}=1-\|\overline{R}\|$, where
$\overline{R}$ is the inverse relation to $R,
$ and if $\overline{\|\cdot\|}$ is not trivial on reflexive
relations, then $\overline{\|\cdot\|}$ is not trivial on
antireflexive relations.}.

For $\alpha,\beta\in S$ we set
$r_{\alpha\beta}=R(\alpha,\beta)+R(\beta,\alpha)$ (thus, $r_{\alpha,\alpha}=2$ for any $\alpha$).

\begin{lemma} If $\bar{x}$ is a positive (i.e $x_i>0$, for $i=\overline{1,n}$)
minimal vector, then $\sum\limits_{i=1}^{n}r_{i\alpha}x_i=\sum\limits_{i=1}^{n}r_{i\beta}x_i$
for any  $\alpha,\beta\in S$.
\end{lemma}

It is clear that  $\sum\limits_{i=1}^{n}r_{i\alpha}x_i=\dis{\p f_R\over\p
x_\alpha}(x_1,\dots,x_n)$. However, to get rid of the condition
$\sum\limits_{i=1}^{n}x_i=1$, we consider the function
$\widehat{f}=\left(\sum\limits_{i=1}^{n}x_i\right)^{-2}f_R$.
It is easy to see that the (unconditional) minimum of $\bar{x}$ is reached on
$\bar{x}$. Hence,
\[
{\p \widehat{f}\over \p x_\alpha}(\bar{x})=
\sum\limits_{i=1}^{n}r_{i\alpha}x_i\left(\sum\limits_{i=1}^{n}x_i\right)^{-2}
-2f_R(\bar{x})\left(\sum\limits_{i=1}^{n}x_i\right)^{-3}=0.
\]
But $\sum\limits_{i=1}^{n}x_i=1$. Therefore
$\sum\limits_{i=1}^{n}r_{i\alpha}x_i=2f_R(\bar{x})$. The right
part does not depend on $\alpha$, whence the statement of lemma
follows.

Elements $\alpha,\beta\in S$ are called {\it twins} if
$s\in S\setminus\{\alpha,\beta\}$ implies $r_{\alpha s}=r_{\beta s}$.

The following statement is almost evident but sometimes useful.

\begin{lemma} If $\bar{x}=(x_1,\dots,x_n)$ is a minimal vector,
$\alpha$ and $\beta$ are twins and $r_{\alpha,\beta}\leq 1$
(certainly, the latter holds if $R$ is a antisymmetrical relation),
then $x_\alpha=x_\beta$.
\end{lemma}

Let $\bar{y}=(y_1,\dots, y_n)$, where $y_i=x_i$ for $i\not\in
\{\alpha,\beta\}$ and $y_\alpha=y_\beta=\dis{1\over 2}(x_\alpha+x_\beta)$.
Then
\[
f_R(\bar{x})-f_R(\bar{y})=x_\alpha^2+x_\beta^2+r_{\alpha,\beta}x_\alpha
x_\beta- {1\over 2}\left(x_\alpha+x_\beta\right)^2-{1\over
4}r_{\alpha,\beta}\left(x_\alpha+x_\beta\right)^2.
\]
We have $\dis{1\over 4}\left(x_\alpha-x_\beta\right)^2$ for $r_{\alpha,\beta}=1,$
and $\dis{1\over 2}\left(x_\alpha-x_\beta\right)^2$ for $r_{\alpha,\beta}=0$.
Consequently, if $x_\alpha\ne x_\beta$, then
$f_R(\bar{x})>f_R(\bar{y})$ that contradicts to minimality $\bar{x}$.

Let $R_1$ and $R_2$ be relations on $S_1$ and $S_2$ respectively,
and let $R_1\sqcup R_2$ be the corresponding relation on $S_1\sqcup S_2,$
the disjoint union of $S_1$ and $S_2$
(if $R_1$ and $R_2$ are relations of partial order, then
($S_1\sqcup S_2$, $R_1\sqcup R_2$) is called the {\it cardinal sum}
$(S_1,R_1)+(S_2,R_2)$ \cite{12}). A set $S$ with a relation $R$
is {\it connected} if $(S,R)\ne (S_1\sqcup S_2,R_1\sqcup R_2)$.

\begin{lemma}{\rm \cite{2}}
\[
\dis \|R_1\sqcup R_2\|^{-1}=\|R_1\|^{-1}+\|R_2\|^{-1}.
\]
\end{lemma}

Let $\bar{z}=(x_1,\dots,x_{n_1},y_1,\dots,y_{n_2})$ be
a minimal vector for $(S_1,R_1)\sqcup (S_2,R_2)$,
$\sum\limits_{i=1}^{n_1}x_i=\lambda$,
$\sum\limits_{i=1}^{n_2}y_i=1-\lambda$.

Then $\bar{x}'=\lambda^{-1}(x_1,\dots,x_{n_1})$
is a minimal vector for $(S_1,R_1)$, and
$\bar{y}'=(1-\lambda)^{-1}(y_1,\dots,y_{n_2})$
is a minimal for $(S_2,R_2)$.
Let $u=\|R_1\|=f_{R_1}(\bar{x}')$, $v=\|R_2\|=f_{R_2}(\bar{y}')$.
Then $\|R\|=\lambda^2u+(1-\lambda)^2v=g(\lambda)$, moreover,
$\lambda$ should yield the minimum $g(\lambda)$ on $[0,1]$.
The derivative with respect to $\lambda$
(consider $u$ and $v$ to be constants)
is
\[
2\lambda u-2(1-\lambda)v=0,\quad \lambda={v\over u+v},\quad
1-\lambda={u\over u+v}.
\]
Substitute this value for $\lambda$ and obtain
\[
f_R(\bar{z})={v^{2}u\over (u+v)^2}+ {u^2v\over (u+v)^2}={uv\over
u+v}=\left({1\over u}+{1\over v}\right)^{-1}.
\]
It is also necessary to check that $\dis{uv\over u+v}<u=g(1)$ and
$\dis{uv\over u+v}<v=g(0)$. For $\lambda=0$ we have
\[
v-{uv\over u+v}={v^2\over u+v}>0,
\]
and, similarly, for $\lambda=1,$ $u^2>0$.
Thus, for
$\lambda=\dis{v\over u+v}$ we really have the
least value for  $\lambda^2u+(1-\lambda)^2v,$ and the lemma is proved.

In connection with this statement, let us introduce the function $\Rho(S)=\Rho(S,R)$
to be equal  $\|R\|^{-1}$.
Then $\Rho(S_1\sqcup S_2, R_1\sqcup R_2)=\Rho(S_1,R_1)+\Rho(S_2,R_2)$.

Let $L_n$ be a linearly ordered set (chain) of the order~$n$.

\begin{lemma} $\dis \Rho(L_n)=1+{n-1\over n+1}.$
\end{lemma}
By lemma 2 $\bar{x} = \dis\left({1\over n},\dots,{1\over n}\right)$
is a minimal vector for $L_n$.
\begin{gather*}
f_R(\bar{x})={1\over n^2}\left(n+C^2_n\right)={1\over
n^2}\left(n+{1\over 2}n(n-1)\right)={n+1\over 2n},\\
\Rho(L_n)={2n\over n+1}=1+{n-1\over n+1}.
\end{gather*}

Thus, $1\leq\Rho<2$ for all linear partially ordered sets.
\begin{remark}\rm
Let $S$ be a not linearly ordered poset.
Call $S$ {\it semilinear} if any element of $S$ is not comparable
with at most one element. It is easy to check that
the following conditions are equivalent:
\begin{enumerate}
\item
$S$ is semilinear;
\item
$S$ is an ordinal sum of antichains of orders 1 and 2~\cite{12};
\item
$\Rho(S)=2$.
\end{enumerate}
\end{remark}

From lemma 3 it immediately follows that
if $R$ is an equivalence relation on $S$,
then $\Rho(S,R)$ is equal to the number of equivalence classes.

The definition of $\Rho$ implies that
if $S'\subset S$ and $R'$ is the relation on $S'$ induced by $R$,
then $\Rho(S',R')\leq\Rho(S,R)$. The set $(S,R)$ is called
{\it $\Rho$-faithful} if $\Rho(S',R')<\Rho(S,R)$ for any
$(S',R')\subset(S,R)$.

It is clear that a complete relation is $\Rho$-faithful only for
$n=1$, and an equivalence relation is $\Rho$-faithful iff each
equivalence class contains one element.

Studying of $\Rho$-faithful sets is reduced to connected sets,
because from lemma 3 it follows that
$(\mathop{\sqcup}\limits_{i=1}^t S_i,\mathop{\sqcup}\limits_{i=1}^t R_i)$
is $\Rho$-faithful iff each $(S_i,R_i)$ is $\Rho$-faithful.

By Lemma 4, $L_n$ is $\Rho$-faithful for any $n$.

The authors do not know a criterion of $\Rho$-faithfulness for
arbitrary binary relations. Such investigation would be seemingly
rather difficult and interesting problem.
In fact, we know a motivation for consideration of the norm of a relation
(and, hence, the function $\Rho(S,R)$ and $\Rho$-faithful sets)
only when $R$ is a partial order. However, the definitions given above
look natural, and it seems that they probably make sense
in general case.

\section{$\pbf{\Rho}$-faithful partially ordered sets}

Here and further $S=(S,\leq)$ is a (finite) partially ordered set
(poset).

Following P.~Gabriel we will say a {\it quiver} instead of
``{\it oriented graph}'' (in general, admitting loops and parallel
arrows, i.e. several arrows between the same vertices). Finite set
$S$ is usually depicted by its Hasse diagram: this is a quiver $Q(S),$
whose two points (elements of $S$) are connected by an arrow $x\longrightarrow y$
if $x<y$ and there is no $z\in S$ such that $x<z<y$.
We usually will draw lines instead of arrows
assuming that an arrow is always directed from the bottom to the top.
The maximal number $w(S)$ of pairwise incomparable elements of $S$
is called the {\it width} of $S$. $S$ is {\it primitive}
if $S=\bigsqcup\limits_{i=1}^tL^i_{n_i}$, i.e.
$S$ is a union of some pairwise incomparable
($x$ is not comparable with $y$ for $x\in L^i$, $y\in L^j$, $i\not=j$) chains.
We also denote a primitive poset by $(n_1,\dots,n_t)$,
in particular $(n)$ is a chain of order $n$.

Representations of posets were introduced in \cite{3} (see
section~5 of present paper). The following criterion of finite
representativity was proved in \cite{7}.

Posets $S$ is finitely represented (i.e. has a finite number of
pairwise nonequivalent indecomposable representations) iff
$S$ does not contain the following subsets: $K_1=(1,1,1,1)$, $K_2=(2,2,2)$,
$K_3=(1,3,3)$, $K_4=(1,2,5)$ and $K_5=(4)\sqcup \widehat{N}$, here
and further we use the notation
$\widehat{N}=$\raisebox{-1.5ex}{$\overset{\overset{\displaystyle
\circ }|}\circ$}$\displaystyle\diagdown$\raisebox{-1.5ex}
{$\overset{\overset{\displaystyle \circ }|} \circ $}.

In \cite{13,8,14,22} infinitely represented posets and other
matrix problems are divided into two types: tame and wild.
In these articles the definitions of tame and wild types formally
are somewhat different, but actually they are coincide,
and we will not remind them here. Note only that
in accordance with \cite{15,22,9} it will be convenient
to differ finitely represented and tame problems.

In \cite{8} it is proved that a poset has tame type iff $S$
contains one of the posets $K_1-K_5$ and does not contain subsets of the following form:
$N_0=(1,1,1,1,1)$, $N_1=(1,1,1,2)$, $N_2=(3,2,2)$,
$N_3=(1,3,4)$, $N_4=(1,2,6)$ and $N_5=(5)\sqcup \widehat{N}$.

\begin{remark}\rm
From two criteria given above it follows that any wild poset
contains a tame subset. It is also true for other matrix
problems but does not follow from the definitions.
\end{remark}

Posets $K_1-K_5$ and $N_0-N_5$ are quite simple
but mysterious in a way (in particular, only $K_5$ and $N_5$ are not primitive).

Below we consider indicated critical posets from the point of view
of the function $\Rho$ defined in section~1.

Using lemmas 4 and 3 it is not difficult to see that
for $i=\overline{1,4},$ $\Rho(K_i)=4$; $\Rho(K_5)$ is also equal to 4,
but for this we need to calculate that $\Rho(\widehat{N})=2,4$. For the sets
$N_0-N_5$, $\Rho(N_i)>4$.

\begin{remark}\rm
Analogously to remark 1 one can verify that the following
conditions are equivalent:
\begin{enumerate}
\item
$S$ does not contain subsets of the form $(1,3)$, $(1,1,1)$, $(2,2)$ and
$\widehat{N}$;
\item
$S$ is an ordinal sum of posets of the form $(1)$, $(1,1)$ and $(1,2)$;
\item
$\Rho(S)< 2,4$.
\end{enumerate}
\end{remark}

In \cite{1,2} it is proved the following statement.

\begin{proposition} Poset $S$ is finitely represented (resp. tame) if
$\Rho(S)<4$ (resp. $\Rho(S)=4$).
\end{proposition}

From the proposition and the finite and tame criteria
it follows that all $\Rho$-faithful sets with $\Rho=4$
are exhausted by the sets $K_1-K_5$. Sets $N_i$ are characterized
in terms of the function $\Rho$ as follows:
$S\in \{N_i\}$ iff
$\Rho(S)>4$, and for any $S'\subset S$, $\Rho(S')\leq 4$.

Remind that every poset of the width $t$ is a union of $t$
disjoint chain (Dilworth's theorem). We call a poset $S$
a {\it wattle} if $S$ is a union of $t$ disjoint chains
$Z_1,\dots, Z_t$, where $|Z_i|\geq 2$, $(i=\overline{1,t}$, $t>1)$,
the minimal element of $Z_i$ is smaller than the maximal element of
$Z_{i+1}$, $(i=\overline{1,t-1})$ and there are no other comparisons
between elements of different chains. Thus, a wattle is given
by the set $\langle n_1,\dots,n_t\rangle$, where $n_i=|Z_i|$,
$(n_i\geq 2)$. The poset $\widehat{N}$ is the simplest wattle
$\langle 2,2\rangle$.

It can be verified that an {\it equal-high} wattle (i.e.
$n_i=\dis{n\over t}$, for $i=\overline{1,t}$) is $\Rho$-faithful, and
the wattle $\{n_1,n_2\}$ with $n_1\ne n_2$ is not $\Rho$-faithful. However,
the wattle $\langle 2,3,2\rangle$ (indicated by M.~Zeldich) is $\Rho$-faithful.

Denote the minimal points of the chains $Z_i$ by $z_i^-$ for
$i=\overline{1,n-1}$ and the maximal points of the chains
$Z_i$ for $i=\overline{2,n}$ by $z_i^+$;
the other points of a wattle are called {\it common points}.

\begin{lemma}
If $\bar{x}=(x(1),\dots,x(n))$ is a positive minimal vector for a wattle $S$,
then there are numbers $\alpha,\beta$ such that
\begin{enumerate}
\item[1)]
$x(s)=\alpha$, if $s$ is a common point;
\item[2)]
$x(z^-_i)+x(z^+_{i+1})=\alpha$ for $i=\overline{1,t-1}$;
\item[3)]
$\sum\limits_{s\in Z_i}x(s)=\beta$ for $i=\overline{1,t}$.
\end{enumerate}
\end{lemma}

A proof immediately follows from lemma 1. (apply the lemma consequently:
firstly to $z_i^-$ and $z_i^+$ $(1<i<t)$ to obtain the property~2),
further to common and uncommon points of every chain $Z_i$ to obtain the property~1),
and at last to common points of different $Z_i$ to obtain the property~3).)

\begin{proposition} If $S$ is a wattle and $|Z_1|=k$, then $S$ can be
$\Rho$-faithful only if $S$ is a uniform wattle, i.e.
\begin{enumerate}
\item[a)]
$k\leq |Z_i|\leq k+1$ for $i=\overline{2,t-1}$.
\item[b)]
$|Z_t|=k$.
\item[c)]
If $m$ is the number of $Z_i'$s such that $|Z_i|=k+1$,
$(0\leq m\leq t-2)$, than $m+1$ and $t$ are relatively primes.
\item[d)]
If $u_1<\dots<u_m\in\{1,\dots,t\}$ are such numbers that
$|Z_{u_i}|\ne k$, then  $u_i=\left[\dis{it\over m+1}\right]+1$.
\end{enumerate}
If conditions $a)$, $b)$, $c)$ and $d)$ are satisfied,
there exists a vector $\bar{x}$ satisfying the conditions of lemma~1.
\end{proposition}

Let $\bar{x}=(x(1),\dots,x(n))$ be a positive minimal vector
for a wattle~$S$. Set $\gamma=x(z_1^-)$. The statements $a)$ and $b)$
immediately follows from lemma~5. It also implies that
$\gamma=\beta-(k-1)\alpha$.
$x(z_2^+)=\alpha-\gamma=k\alpha-\beta$.
If $t>2$, then
\[
\sum\limits_{s\in
Z_2}x(s)=\beta=r_2\alpha+\alpha-\gamma+x(z_2^-),\]
where $r_2=k-2$ for $|Z_2|=k$
and $r_2=k-1$ for $|Z_2|=k+1$.

Thus, for $|Z_2|=k$
\[
x(z_2^-)=\beta-(k-2)\alpha -\alpha + \gamma=2\gamma,
\]
and for $|Z_2|=k+1$
\[ x(z_2^-)=2\gamma -\alpha.\]

Analogously, we obtain that for $i=\overline{2,t-1}$
\[x(z_i^-)=i\gamma-u_i\alpha,\] where $u_i$ is the number of such
$j=\overline{2,i}$ that $|Z_j|=k+1$. From $0<x(z_i^-)<\alpha$
it follows that $u_i=\left[\dis{i\gamma\over\alpha}\right]$.

On the other hand,
\begin{gather*}
(k-1)\alpha+\gamma=\beta=\sum\limits_{s\in
Z_i}x(s)=(k-1)\alpha+x(z_n^+)=\\ =(k-1)\alpha+\alpha-x(z_{t-1}^-)=
k\alpha-(t-1)\gamma+m\alpha, \\ (m=u_{t-1}),\quad
t\gamma=\alpha(1+m),\quad \gamma={(1+m)\alpha\over t},\quad
u_i=\left[i(1+m)\over t\right],
\end{gather*}
that implies the statements $c$ and $d$ (if $c$ would not hold,
then for some $i$ $x(z_i^-)=0$). If the conditions $a)$, $ b)$, $c)$
and $d)$ are satisfied, from the above it follows that a positive
vector  $\bar{x}$ satisfying the conditions of lemma~6 can be uniquely
constructed.

The proposition implies that for given $t$, $k$, $m$ and
a wattle given by the conditions $a$, $b$, $d$,
there exists a positive vector $\overline{x}$
satisfying the conditions of lemmas 5 and 1
iff condition $c$ is satisfied.

\begin{example}\rm
Let $t=7$, $m=3$ and $k=2$.
Then $u_1=2$, $u_2=4$, $u_3=6$, $Z=\langle 2,3,2,3,2,3,2\rangle$.
Thus, the condition $d$ guaranties a uniform distribution
of long chains among all chains $Z_i$.
\end{example}

\medskip

\noindent {\bf Conjecture 1.} {\it Every $\Rho$-faithful connected
partially ordered set is either a chain, or a uniform wattle.}

\medskip

Using conjecture 1 and proposition 1 it is not difficult to
reestablish the finiteness and tameness criteria, i.e. lists
$K_i$ and $N_i$.

{\it Remark to the English translation.} A proof of Conjecture~1
is published in the article of A.~Sapelkin ``$\Rho$-faithful
partially ordered sets'' (Ukr. Math. J., 2002, N~7).

\section{Numeral function $\pbf{\rho}$ \\
and a characterization of Dynkin schemes}

Denote by  ${\mathbb N}^\infty$ a set consisting of natural numbers and
symbol $\infty$ and let $\infty+n=\infty-n
=\infty+\infty=\infty$ $(n\in {\mathbb N})$.
We construct a function $\rho$ on ${\mathbb N}^\infty$.
Identifying natural number $n$ with the chain $L_n$
we set (taking into account lemma 4)
$\rho(n)=1+\frac{n-1}{n+1}$ and $\rho(\infty)=\lim\limits_{n\to\infty}\rho(n)=2$.
Sometimes we also use (according to the same formula) $\rho(0)=0$.

We assume $\rho(n_1,\ldots,n_t)=\sum\limits_{i=1}^t\rho(n_i)$
$(n_i\in {\mathbb N}^\infty)$.

We shall often consider the condition $\rho(n_1,\ldots,n_t)\leq 4$.
For $t\leq 2$ it always holds, and the equality takes place only
for $n_1=n_2=\infty$; for $t\geq 4$ the condition (with equality)
is satisfied only for $t=4$ and $n_1=n_2=n_3=n_4=1$.
For $t=3$, if $n_1=\infty$
(assume $n_1\geq n_2\geq n_3$), only the case
$n_2=n_3=1$ is possible. Thus, the condition is not evident
only for $t=3$ and $n_i<\infty$.

The following lemma follows from a trivial algebraic transformation.

\begin{lemma} For $n_i<\infty$ the following conditions are equivalent:
\begin{enumerate}
\item[i)] $\rho(n_1,n_2,n_3)<4$ (resp. $\rho(n_1,n_2,n_3)=4$);

\item[ii)] $(n_1+1)^{-1}+(n_2+1)^{-1}+(n_3+1)^{-1}>1$ (resp.
$(n_1+1)^{-1}+(n_2+1)^{-1}+(n_3+1)^{-1}=1)$;

\item[iii)] $\mu(n_1-1,n_2-1,n_3-1)<4$ (resp.  $\mu(n_1-1,n_2-1,n_3-1)=4$).
\end{enumerate}

The equation $\rho(x_1,\ldots,x_t)=4$ $(x_i \in{\mathbb N}^\infty)$
has exactly six solutions: $(\infty,\infty)$, $(\infty,1,1)$,
$(5,2,1)$, $(3,3,1)$, $(2,2,2)$ and $(1,1,1,1)$.
\end{lemma}

From the equivalence of i) and ii) it follows that well known
statements using $\sum\limits_{i=1}^3 n_i^{-1}$ can be reformulated
in terms of $\rho$. For example, if $G$ is a group with generators
$g_1$, $g_2$, $g_3$ and relations
$g_1^2=g_2^2=g_3^2=(g_1g_2)^{n_1}= (g_1g_3)^{n_2}=(g_2g_3)^{n_3}$
(see \cite{6}), then $G$ is finite iff $\rho(n_1-1,n_2-1,n_3-1)<4$,
and in this case $|G|=8(4-\rho(n_1-1,n_2-1,n_3-1)^{-1}$. Of course,
there is no sense in such reformulations but they confirm
the prevalence of~$\rho$.

For $t=3$ and $n_i<\infty$ the most suitable function is $\mu$,
however, although a transfer of $\mu$ and $\sum n_i^{-1}$
to general case is possible, it does not seem to be natural.

Let $\Gamma$ be a graph. We will consider only finite graphs
but admit (unlike, for example,~\cite{6}) loops and parallel
edges (i.e. some edges between the same vertices).
Thus, $\Gamma$ consists of two (finite) sets
$\Gamma_v$ ($\not=\varnothing$) and $\Gamma_e$ and
a map $\varphi$ associating one-elements (if $r$ is a loop)
or two-elements subset of a set of vertices $\Gamma_v$
for every edge $r\in \Gamma_e$.
By $\Gamma_L$ denote the set of loops of $\Gamma_e$.
Every graph is a union of its {\it connected components}.

For $x\in \Gamma_v$ we set $\varphi^{-1}(x)=\{\alpha\in\Gamma_e
\;|\;\varphi(\alpha)\ni x\}$. The {\it degree} $g(x,\Gamma)=g(x)$
of vertex $x$ is $|\varphi^{-1}(x)|$,
$x$ is a {\it branch point} if $g(x)>2$.
By $A_l$, $(l\geq 1)$, denote the graph
$\begin{array}{@{}c@{\!\!}c@{}c@{}c@{\!}c@{}}
\overset{a_1}{\circ}&\mbox{---}&\,\cdots
&\mbox{---}&\overset{a_l}{\circ}
\end{array}$ ($g(a_1)=g(a_l)=1$, $g(a_i)=2$ for $1<i<l$).
We call this graph {\it simple}.

In essence, the connection of ``finite'' solutions of the equation
$\rho(x_1,\ldots,x_t)=4$ and Dynkin schemes and extended schemes
is well known (in terms ii). If $\Gamma$ is a tree
(i.e. connected graph without cycles)~\cite{6}
with unique branch point $x$ and
$\Gamma_1,\ldots,\Gamma_t$ $(t=g(x)\geq 3)$
are the connected components of $\Gamma'=\Gamma\backslash x$
($\Gamma'_v=\Gamma_v\backslash x$, $\Gamma'_e=\Gamma_e\backslash
\varphi^{-1}(x)$), then $\Gamma$ is an extended Dynkin scheme
($\widetilde{E_8}$, $\widetilde{E_7}$,
$\widetilde{E_6}$ or $\widetilde{D_4}$~--- see list II below)
if and only if $\rho(|\Gamma_1|,\ldots,|\Gamma_t|)=4$.
In this section we characterize all Dynkin schemes and
extended schemes in terms of the function $\rho$.

For $\alpha \in \Gamma_e$ we denote the graph
$\Gamma\backslash \alpha$ ($(\Gamma'_\alpha)_v=\Gamma_v$,
$(\Gamma'_\alpha)_e=\Gamma_e\backslash \{\alpha\}$) by $\Gamma'_\alpha$.
An edge $\alpha$ of a connected graph $\Gamma$ is {\it cyclic} if
$\Gamma'_\alpha$ is connected. Otherwise, the graph $\Gamma'_\alpha$
falls into two connected components $\Gamma'(x,\alpha)\ni x$ and
$\Gamma'(y,\alpha)\ni y$, where  $\varphi(\alpha)=\{x,y\}$.

$(v,f)$-{\it graph} $\overline{\Gamma}=(\Gamma,v,f)$ is
a graph $\Gamma$ together with functions $v$ and $f$ given on the sets
$\Gamma_v$ and $\Gamma_e$ and taking the values in ${\mathbb N}^\infty$
and $R[1,\infty]$ respectively, where $R[1,\infty]$
consists of real numbers being more or equal 1 and symbol $\infty$.

If the values of the function $v$ (resp. $f$) for all vertices
(resp. edges) are equal~1, then $(\Gamma,f)$ is an {\it $f$-graph}
(resp. $(\Gamma,v)$ is a $v$-graph)\footnote{In fact, we consider
only $f$-graphs (in sections~3, 4) and $v$-graphs (in~section~8),
but it seems natural to give some definitions for
$(v,f)$-graphs.}. If $f(\alpha)\in {\mathbb N}^\infty$ $(\alpha\in
\Gamma_e)$, then $(\Gamma,f)$ is an {\it integral} $f$-graph; if,
additionally, $f(\alpha)\geq 3$ and $\Gamma$ does not contain
loops and parallel edges, then $(\Gamma,f)$ is called a {\it
Coxeter graph} \cite{6}.

{\it $f$-degree} $g_f(x)$ of a vertex $x$ of a $(v,f)$-graph
$(\Gamma,v,f)$ is $\sum\limits_{\alpha\in \varphi^{-1}(x)} f(\alpha)$.

An {\it incident pair} $(x,\alpha)$
($x\in\Gamma_v$, $\alpha\in\Gamma_e$, $x\in\varphi(\alpha)$)
is {\it simple} if $\alpha$ is not cyclic, $g(y)\leq 2$ ($\varphi(\alpha)=\{x,y\}$)
and the graph $\Gamma'(y,\alpha)$ is simple;
a pair $(x,\alpha)$ is {\it $(v,f)$-simple}
(in $\overline{\Gamma}=(\Gamma,v,f)$), if, additionally,
$v(a)=1$ for $a\in \Gamma'(y,\alpha)$, $g(a,\Gamma)=g(a)\not=1$,
and, finally, $f(\beta)=1$ for all $\beta\in \Gamma'(y,\alpha)_e$.

We assign
the number $\displaystyle \frac{\partial x}{\partial \alpha}$
to a pair  $(x,\alpha)$ ($\varphi(\alpha)\ni x$)
in the following way:

1) $\displaystyle \frac{\partial x}{\partial \alpha}=
\rho\left(\sum\limits_{z\in \Gamma'(y,\alpha)} v(z)\right)$,
if $(x,\alpha)$ is a $(v,f)$-simple pair;

2) $\displaystyle \frac{\partial x}{\partial \alpha}=2$
(i.e. $\rho(\infty)$), if $\alpha\not\in\Gamma_L$ and
$(x,\alpha)$ is not a $(v,f)$-simple pair;

3) $\displaystyle \frac{\partial x}{\partial \alpha}=4$,
if $\alpha\in \Gamma_L$.

If $x$ is a vertex of a $(v,f)$-graph $\overline{\Gamma}$,
then its {\it $\rho$-degree} $\displaystyle
g_\rho(x)=\rho(v(x)-1)+\sum\limits_{\alpha\in\varphi^{-1}(x)}
f(\alpha) \frac{\partial x}{\partial \alpha}$.

If $\overline{\Gamma}=(\Gamma,f)$ is an $f$-graph, (or $v(x)=1$),
then the summand $\rho(v(x)-1)$ vanishes.

Of course, the definition of a $\rho$-degree also makes sense
when $v$ and $f$ are identically equal~1.

\begin{remark}\rm
For a $v$-graph $\overline{\Gamma}=(\Gamma,v)$ it is possible to
construct a graph $G$ in the following way. For each $x\in \Gamma_v$
such that $1<v(x)<\infty$ we add points $a_2^x,\ldots,a_m^x$,
where $m=f(x)$  and
$x\mbox{---}a_2^x\mbox{---}\cdots\mbox{---}a_m^x$ is a simple
graph of $G$ (and there are no other edges,
which are incident to vertices $a_2^x,\ldots,a_m^x$ in $G$),
and for each $y$, if $f(y)=\infty$,
we add points  $b_1^y$ and $b_2^y$ such that
$b_1^y\mbox{---}y\mbox{---} b_2^y$ (and there are no other edges,
which are incident to $b_1^y$ and $b_2^y$ in $G$) (see~\cite{9}).
The definition of $\rho$-degree implies that if $z\in \Gamma_v$,
then its $\rho$-degree in $v$-graph $(\Gamma,v)$ and in graph $G$
coincide, and if $w\in G_v\backslash \Gamma_v$, then $g_\rho(w)< 4$.
\end{remark}

\begin{example}\rm
 $\Gamma= \underset{x}{\circ}\!\overset{\alpha}{\mbox{---}}\underset{y}{\circ}$;
$v(x)=3$, $v(y)=\infty$ in $v$-graph $\overline{\Gamma}$;
$\displaystyle g_\rho(x)=\rho(3-1)+\frac{\partial x}{\partial
\alpha}=\rho(2)+\rho(\infty)=3\frac 13$, $\displaystyle
g_\rho(y)=\rho(\infty-1)+\frac{\partial y}{\partial
\alpha}=\rho(\infty)+\rho(3)=3\frac 12$.

$G=\underset{a_3^x}{\circ}\!\!\mbox{---}\!\!\underset{a_2^x}{\circ}\!\!\!\mbox{---}\!
\underset{x}{\circ}\!\!\overset{\alpha}{\mbox{---}}\!\!\!\!\!\!
\overset{\overset{\mbox{\scriptsize\raisebox{-0.9mm}[0pt][0pt]{$\;
\; \;\,
 {\displaystyle \circ}{ b_1^y}$}}}{|}}{\underset{y}{\circ}}
\!\!\!\!\!\!\mbox{---}\!\underset{b_2^y}{\circ}$ $\displaystyle
\frac{\partial x}{\partial \alpha}=2$ (the pair $(x,\alpha)$ is not
simple since though $G'(y,\alpha)\simeq A_3$, but $g(y)> 2$),
$\displaystyle \frac{\partial x}{\partial
(x,a_2^x)}=\rho(2)=1\frac 13$, $\displaystyle g_\rho(x)=3\frac
13$, $\displaystyle \frac{\partial y}{\partial
\alpha}=\rho(3)=1\frac 12$, $\displaystyle \frac{\partial
y}{\partial (y,b_1^y)}=\frac{\partial y}{\partial (y,b_2^y)}=1$,
$\displaystyle g_\rho(y)=3\frac 12$, $\displaystyle
g_\rho(a_2^x)=\frac{\partial a_2^x}{\partial (a_2^x,a_3^x)}+
\frac{\partial a_2^x}{\partial (a_2^x,x)}=1+2=3$,
$g_\rho(b_1^y)=g_\rho(b_2^y)=g_\rho(a_3^x)=2$.
\end{example}

Below we give the list of Dynkin schemes

\medskip

$A_l$ \ \
$\begin{array}{c}
\mbox{\includegraphics{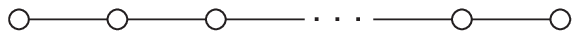}}
\end{array}$ \ \ ($l\geq 1$ vertices)

\medskip

$B_l$ \ \
$\begin{array}{c}
\mbox{\includegraphics{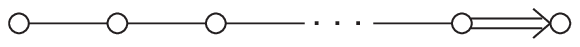}}
\end{array}$ \ \ ($l\geq 2$ vertices)

\medskip

$C_l$ \ \
$\begin{array}{c}
\mbox{\includegraphics{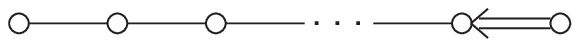}}
\end{array}$ \ \ ($l\geq 3$ vertices)

\medskip

$D_l$ \ \
$\begin{array}{c}
\mbox{\includegraphics{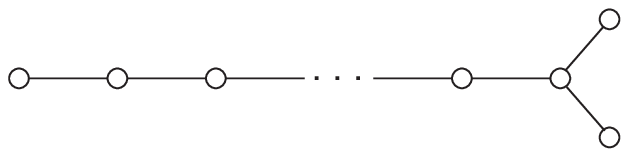}}
\end{array}$ \ \ ($l\geq 4$ vertices)

\medskip

$E_6$ \ \ $\begin{array}{c}
\mbox{\includegraphics{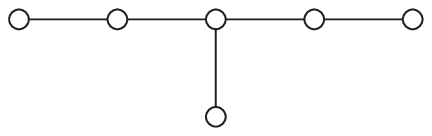}}
\end{array}$\hfill (I)

\medskip

$E_7$ \ \
$\begin{array}{c}
\mbox{\includegraphics{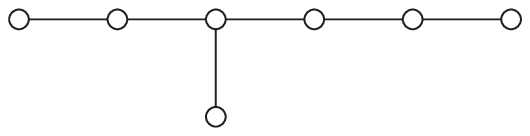}}
\end{array}$

\medskip

$E_8$ \ \
$\begin{array}{c}
\mbox{\includegraphics{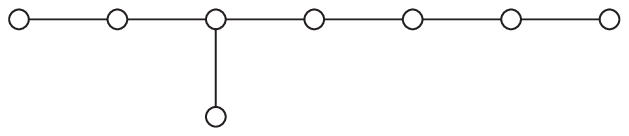}}
\end{array}$

\medskip

$F_4$ \ \ $\begin{array}{c}
\mbox{\includegraphics{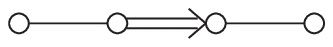}}
\end{array}$

\medskip

$G_2$ \ \ $\begin{array}{c}
\mbox{\includegraphics{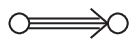}}
\end{array}$

\medskip

\noindent
and extended Dynkin schemes

\medskip

$\widetilde{A_{l-1}}$ \ \
$\begin{array}{c}
\mbox{\includegraphics{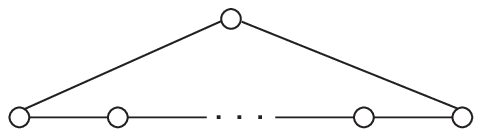}}
\end{array}$ \ \ ($l\geq 1$ vertices)

\medskip

$\widetilde{C_{l}}$ \ \
$\begin{array}{c}
\mbox{\includegraphics{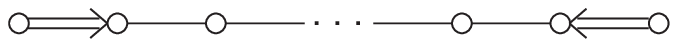}}
\end{array}$ \ \ ($l+1$ vertices, $l\geq 2$)

\medskip

$\widetilde{D_{l}}$ \ \
$\begin{array}{c}
\mbox{\includegraphics{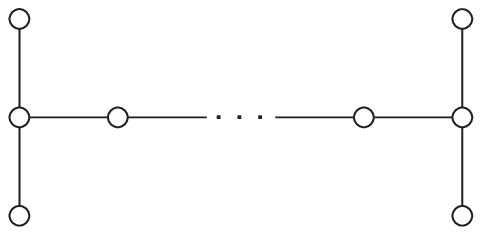}}
\end{array}$ \ \ ($l+1$ vertices, $l\geq 4$)

\medskip

$\widetilde{B_{l}}$ \ \
$\begin{array}{c}
\mbox{\includegraphics{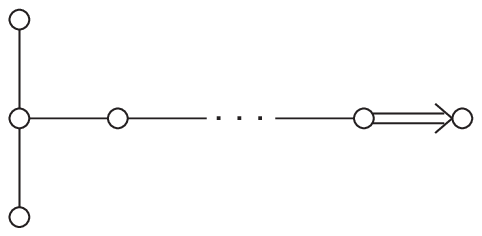}}
\end{array}$\ \ ($l+1$ vertices, $l\geq 3$)\hfill (II)

\medskip

$\widetilde{G_{2}}$ \ \ $\begin{array}{c}
\mbox{\includegraphics{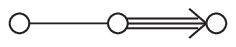}}
\end{array}$

\medskip

$\widetilde{E_{6}}$ \ \ $\begin{array}{c}
\mbox{\includegraphics{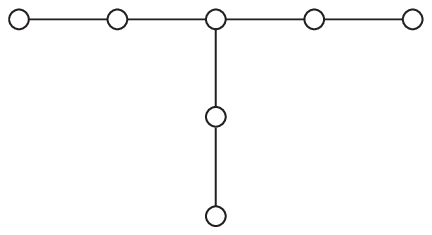}}
\end{array}$

\medskip

$\widetilde{E_{7}}$ \ \
$\begin{array}{c}
\mbox{\includegraphics{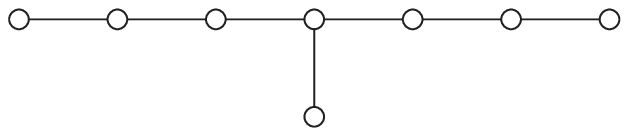}}
\end{array}$

\medskip

$\widetilde{E_{8}}$ \ \
$\begin{array}{c}
\mbox{\includegraphics{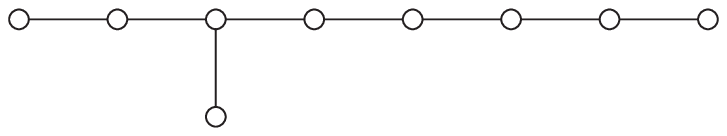}}
\end{array}$

\medskip

$\widetilde{F_{4}}$ \ \
$\begin{array}{c}
\mbox{\includegraphics{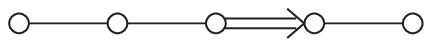}}
\end{array}$

\medskip

$BA_2$ \ \ $\begin{array}{c}
\mbox{\includegraphics{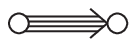}}
\end{array}$

\medskip

$FE_{6}$ \ \ $\begin{array}{c}
\mbox{\includegraphics{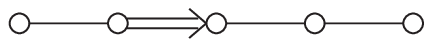}}
\end{array}$

\medskip

$BA_{l+1}$ \ \ $\begin{array}{c}
\mbox{\includegraphics{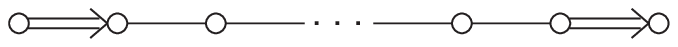}}
\end{array}$ \ \ ($l+1$ vertices, $l\geq 2$)

\medskip

$BD_{l}$ \ \
$\begin{array}{c}
\mbox{\includegraphics{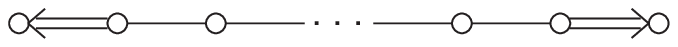}}
\end{array}$ \ \ ($l+1$ vertices, $l\geq 2$)

\medskip

$CA_{l}$ \ \
$\begin{array}{c}
\mbox{\includegraphics{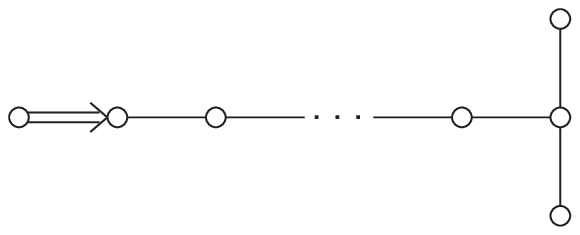}}
\end{array}$ \ \ ($l+1$ vertices, $l\geq 3$)

\medskip

$GD_4$ \ \ $\begin{array}{c}
\mbox{\includegraphics{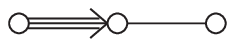}}
\end{array}$

\medskip

Each Dynkin scheme (or extended scheme)
can be considered as an $f$-graph,
in which for each edge $\alpha$ an additional
arrow is drawn if $f(\alpha)\ne 1$ (then $f(\alpha)$
is equal to the multiplicity of the corresponding arrow).

An orientation of an arrow is essential for a scheme
(reorientation changes the system of roots, Lie algebras and so on).
However, inspecting the lists I and II it easy to see that a scheme
is remaining in the list after reorientation. We will say that
an $f$-graph $G$ {\it generates} a Dynkin scheme (resp. an extended Dynkin scheme)
if we obtain a Dynkin scheme (resp. an extended Dynkin scheme)
for some (and, therefore, for any) arrangement of arrows on edges
$\alpha$ such that $f(\alpha)>1$.

We sometimes will write $f(\alpha)\ne 1$ over the edge $\alpha$
or near it on the picture.

The $f$-graph corresponding to a scheme $X$
we denote by $\underline{X}$
(for example, $\underline{BA_2}$ is $\begin{array}{@{}c@{}c@{}c@{}}
{\circ}&\overset{4}{{\mbox{---}}} &{\circ}
\end{array}$).

$\underline{\rm I}=\{A_l,\underline{B_l},D_l,E_6,E_7,E_8,
\underline{F_4},\underline{G_2}\}$ ($\underline{C_l}\sim
\underline{B_l}$) corresponds to the list I, and
$\underline{\rm II}=\{\widetilde{A_{l-1}},\widetilde{\underline{C_l}},
\widetilde{\underline{D_l}}, \widetilde{\underline{B_l}},
\widetilde{\underline{G_2}}, \widetilde{E_6},\widetilde{E_7},
\widetilde{E_8},\widetilde{\underline{F_4}}, \underline{BA_2}\}$
to the list II.

\begin{proposition}
A connected integral $f$-graph $\overline{\Gamma}$ generates a Dynkin scheme
iff $\rho$-degree of any its vertex is less than 4;
$\overline{\Gamma}$ generates an extended Dynkin scheme iff
$\max g_\rho(x)=4$ $(x\in\Gamma_v)$.
\end{proposition}

Firstly, suppose that $\Gamma$ contains a loop.
In the lists I, II there is
the only such graph $\widetilde{A_0}$: $\Gamma_v=\{x\}$,
$\Gamma_e=\{\alpha\}$, $\varphi(\alpha)=x$ $(f(\alpha)=1)$.
By our formulas $\displaystyle\frac{\partial x}{\partial
\alpha}=4$ $g_\rho(x)=4$, i.e.
$\widetilde{A_0}$ satisfies the condition of proposition~2.

Let $\alpha$ be a loop at a vertex $x$ and
$\overline{\Gamma} \not \sim\widetilde{A_0}$.
Then either $x$ is connected with a vertex
$y\not= x$, or there is one more loop at $x$, or $f(\alpha)\not=1$.
In any case ${g_\rho}(x)>4$.

Therefore, further we may assume $\Gamma_L=\varnothing$.

Suppose that $\Gamma$ does not contain a cycle $\begin{array}{c}
\mbox{\includegraphics{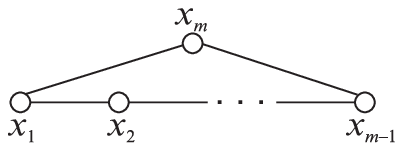}}
\end{array}$ $(m>1)$.
There is the only such graph $\widetilde{A_{l-1}}$, $l> 1$ in the lists
$\underline{\rm I}$, $\underline{\rm II}$. For this graph, for any
$x\in \Gamma_v$ $|\varphi^{-1}(x)|=2$, and both edges
$\alpha,\beta\in \varphi^{-1}(x)$ are cyclic. $\displaystyle
\frac{\partial x}{\partial \alpha}=2$, $\displaystyle
\frac{\partial x}{\partial \beta}=2$, $\displaystyle {g_\rho}(x)=
\frac{\partial x}{\partial \alpha}+\frac{\partial x}{\partial
\beta}=4$. If $\overline{\Gamma}$ is an ``other'' (connected)
integral $f$-graph containing a cycle, then this cycle contains
such point $x$, that either $g(x)>2$, or $\varphi(\alpha)\ni x$
and $f(\alpha)\not=1$. In both cases ${g_\rho}(x)>4$.

Thus, we may assume that graph $\Gamma$ is acyclic
(i.e. ``a tree'').

Let $g_f(\overline{\Gamma})=\max\limits_{x\in \Gamma_v} g_f (x)$.
For $f$-graph corresponding $\underline{\rm I}$, $\underline{\rm
II}$, $g_f(\overline{\Gamma})\leq 4$. Clearly, if $g_f(x)>4$, then
${g_\rho}(x)>4$, i.e. such $f$-graph does not satisfy the
conditions of the proposition.

1) Let $g_f(\overline{\Gamma})=4$.
In $\underline{\rm I}$, $\underline{\rm II}$,
there are the following $f$-graphs corresponding to this case:

\begin{center}
$\begin{array}{c}
\mbox{\includegraphics{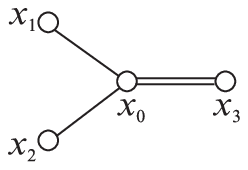}}
\end{array}$ \ \ $(\widetilde{\underline{B_3}})$
\qquad
$\begin{array}{c}
\mbox{\includegraphics{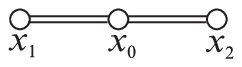}}
\end{array}$ \ \ $(\widetilde{\underline{C_2}})$\\
$\begin{array}{c}
\mbox{\includegraphics{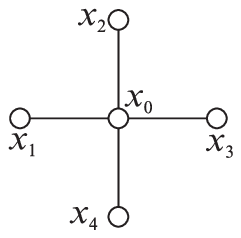}}
\end{array}$ \ \ $(\widetilde{D_4})$\qquad
$\begin{array}{c}
\mbox{\includegraphics{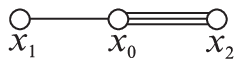}}
\end{array}$ \ \ $(\widetilde{\underline{G_2}})$ \\
$\begin{array}{c}
\mbox{\includegraphics{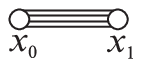}}
\end{array}$ \ \ $(\underline{BA_2})$
\end{center}

Everywhere $g_f(x_0)=4$ (and, moreover, $g_f(x_1)=4$ in
$\underline{BA_2}$). In these cases $\displaystyle\frac{\partial
x_0}{\partial \alpha}=1$, ${g_\rho}(x_0)= g_f(x_0)=4$
($\varphi(a)\ni x_0$). If $g_f(x)=1$, then ${g_\rho}(x)\leq 2$
($x_1$, $x_2$ in $\widetilde{\underline{B_3}}$; $x_1$, $x_2$,
$x_3$, $x_4$ in $\widetilde{D_4}$; $x_1$ in
$\widetilde{\underline{G_2}}$). If $g_f(x)=2$, then
${g_\rho}(x)\leq 4$ ($x_3$ in $\widetilde{\underline{B_3}}$;
$x_1$, $x_2$ in $\widetilde{\underline{C_2}}$). In
$\widetilde{\underline{G_2}}$ $g_\rho(x_2)=3\rho(2)$,
$\rho(2)=1\frac 13$, $g_\rho(x_2)=4$.
Thus, the listed above $f$-graphs
satisfy the conditions of the proposition.

Now let $\overline{\Gamma}$ be an arbitrary $f$-graph
satisfying the conditions of the proposition; $\Gamma_v\ni x$,
$g_f(x)=4$. Then ${g_\rho}(x)=4$, and $\displaystyle\frac{\partial
x}{\partial \alpha}=1$
for any $\alpha\in \varphi^{-1}(x)$
$\displaystyle
\Bigg({g_\rho}(x)=\sum\limits_{\alpha\in\varphi^{-1}(x)} f(\alpha)
\frac{\partial x}{\partial \alpha}
=g_f(x)=\sum\limits_{\alpha\in\varphi^{-1}(x)} f(\alpha), \ 1\leq
\frac{\partial x}{\partial \alpha}\Bigg)$.
We obtain $\widetilde{D_4}$ for $g(x)=4$,
$\widetilde{\underline{B_3}}$ for $g(x)=3$,
either $\widetilde{\underline{C_2}}$, or $\widetilde{\underline{G_2}}$
for $g(x)=2$, and, finally, $\underline{BA_2}$ for $g(x)=1$.

2) Let $g_f(\overline{\Gamma})=3$ and let $|\{x\in
\Gamma_v|g_f(x)=3\}|=1$. In $\underline{\rm I}$, $\underline{\rm
II}$ it is one of $\underline{B_l}$ for  $l>2$, $D_l$, $E_6$,
$E_7$, $E_8$, $\widetilde{E_6}$, $\widetilde{E_7}$ or
$\widetilde{E_8}$. It is easy to see that for $y\not= x$ $(g_f
(y)\leq 2)$ ${g_\rho}(y)< 4$.
\begin{gather*}
\underline{B_l}: \ {g_\rho}(x)=2+1+\frac{l-3}{l-1}<4; \\ D_l: \
{g_\rho}(x)=1+1+1+\frac{l-4}{l-2}<4;\\ E_6: \ {g_\rho}(x)=1+1+1+0
+\frac 13 +\frac 13=3\frac 23;
\\
 E_7: \
{g_\rho}(x)=1+1+1+0+\frac13+\frac 12=3\frac 56;
\\
 E_8: \
{g_\rho}(x)=1+1+1+0 +\frac 13 +\frac 35=3\frac{14}{15};
\\
\widetilde{E_6}: \ {g_\rho}(x)=1+1+1+\frac 13+\frac13+\frac 13=4;
\\
\widetilde{E_7}: \ {g_\rho}(x)=1+1+1+0 +\frac 12 +\frac 12=4; \\
\widetilde{E_8}: \ {g_\rho}(x)=1+1+1+0 +\frac 13+\frac23=4.
\end{gather*}

If the conditions of the proposition hold for $\overline{\Gamma}$
and $\varphi(\alpha)\ni x$, then $f(\alpha)\leq 2$.
If $f(\alpha)=2$, then $\varphi(\alpha)=\{x,y\}$ and $g(y)=1$
(otherwise $g_f(y)=3$),
$\Gamma\backslash\{x,y\}=A_m$
(otherwise there is $z\in \Gamma_v\backslash x$ such that $g_f(z)\geq 3$),
moreover, if $\varphi(\beta)=\{x,z\}$, then the pair $(x,\beta)$ is $(v,f)$-simple.
We have $\underline{B_l}$. If for all $\alpha\in \Gamma_e$,
$f(\alpha)=1$, then $\Gamma\backslash x= A_{n_1}\amalg
A_{n_2}\amalg A_{n_3}$, $\displaystyle
g_\rho(x)=3+\frac{n_1-1}{n_1+1}+\frac{n_2-1}{n_2+1}+
\frac{n_3-1}{n_3+1}$ $(n_i \in{\mathbb N})$.
From the condition ${g_\rho}(x)\leq4$ we obtain $E_6$, $E_7$, $E_8$,
$\widetilde{E_6}$, $\widetilde{E_7}$, $\widetilde{E_8}$ (see lemma~6).

3) Let $g_f(\overline{\Gamma})=3$
$|\{x\in\Gamma_v|g_f(x)=3\}|=2$, $g_f(x)=g_f(y)=3$, and
let $\overline \Gamma$ satisfy the conditions of the proposition.
Suppose $g(x)=3$. $\Gamma\backslash x$ has 3 components
$\Gamma_1$, $\Gamma_2$, $\Gamma_3$,  and one of them contains $y$.
Therefore, $\displaystyle\frac{\partial x}{\partial
\alpha_1}=2$, $\displaystyle\frac{\partial x}{\partial
\alpha_2}=1$, $\displaystyle\frac{\partial x}{\partial
\alpha_3}=1$ $(\varphi(\alpha_i)\ni x)$.

If $g(y)$ is also 3, then we get $\widetilde{D_l}$ for $l\geq 5$
(taking into account that $g_f(z)\leq 2$ for $z\not\in x,y$).
If $g(y)=2$, then we have $\widetilde{\underline{B_l}}$ for $l\geq 4$.

If $g(x)=g(y)=2$, then two cases are possible:

a) $\varphi(\alpha)=\{x,y\}$, $f(\alpha)=2$;

b) $\varphi(\alpha)=\{x,u\}$, $\varphi(\beta)=\{y,v\}$,
$f(\alpha)=f(\beta)=2$ (and $f(\gamma)=1$ for
$\gamma\not\in\{\alpha,\beta\}$).

In both cases $\Gamma=A_m$ (otherwise, there is $z\not\in \{x,y\}$ and $g_f(z)=3$).

a) $x=a_t$, $y=a_{t+1}$, $t>1$, $t+1<m$. Then $\displaystyle
\frac{\partial x}{\partial \alpha_{t-1,t}}=\rho(t-1)$ (where
$\alpha_{t-1,t}\in \Gamma_e$,
$\varphi(\alpha_{t-1,t}=\{a_{t-1},a_t\}$), $\displaystyle
\frac{\partial x}{\partial \alpha}=\rho(m-t)$,
${g_\rho}(x)=\rho(t-1)+2\rho(m-t)\leq 4$.

By lemma~6 ${g_\rho}(x)\leq 4$ iff $\mu(t-2,m-t-1, m-t-1)\leq 4$.
If $t-2\not=0$, $m-t-1\not=0$, then $t-2=m-t-1=1$, $t=3$, $m=5$.
$\overline{\Gamma}=
\begin{array}{@{}c@{}c@{}c@{}c@{}c@{}c@{}c@{}c@{}c@{}}
{\circ}&\mbox{---}&
{\circ}&\mbox{---}&{\circ}&\overset{2}{\mbox{---}}&{\circ}&\mbox{---}&{\circ}
\end{array}=\widetilde{\underline{F_4}}$.
If  $t-2=0$, then $m-t-1\leq 2$, and we have $t=2$, $m\in
\{4,5\}$, i.e. either $\underline{F_4}$ or $\underline{\widetilde{F_4}}$.
For $m-t-1=0$ we replace $x$ by $y$.

b) $x=a_2$, $y=a_{m-1}$ (otherwise, there is $z$ such that ${g_\rho}(z)=3$),
and we have $\widetilde{\underline{C_l}}$, $l>2$.

Let $g(x)=1$, then $\Gamma_v =\{x,y\}$, and we have
$\underline{G_2}$.

Thus, we have obtained all schemes of the case~3. Clearly,
${g_\rho}(y)\leq 4$ and $g_\rho(z)\leq 4$ for $z\not\in(x,y)$.

4) $g_f(\overline{\Gamma})=3$, $|\{x\in\Gamma_v|g_f(x)=3\}|>2$.

There are no such schemes in $\underline{I}$, $\underline{II}$.
Let $g_f(x)=g_f(y)=g_f(z)=3$, and suppose there is a path $x\cdots
\overset{\alpha}{\mbox{---}}y\overset{\beta}{\mbox{---}}\cdots z$.
Then $\displaystyle\frac{\partial y}{\partial
\alpha}=\frac{\partial y}{\partial \beta}=2$ and since
$g_f(y)=3$, there is one more edge $\gamma$, $\varphi(\gamma)\ni
y$, $f(\alpha)\not=1$, or $f(\beta)\not=1$. In any case ${g_\rho}(y)>4$.

5) $g_f(\overline{\Gamma})\leq 2$. We have $A_l$,
$\underline{B_2}$.

We have considered all schemes of  $\underline{I}$,
$\underline{II}$ in our proof. Hence, the proposition~2 is proved.

\section{Characterization of Coxeter graph}

In \cite{6} it is shown that if (and only if) $(W,S)$ is an
irreducible finite system of Coxeter, then its Coxeter graph is
isomorphic to one of $f$-graphs of the list III:

\medskip

$A_l$ \ \
$\begin{array}{c}
\mbox{\includegraphics{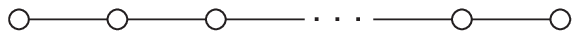}}
\end{array}$ \ \ ($l\geq 1$ vertices)

\medskip

$B_l$ \ \
$\begin{array}{c}
\mbox{\includegraphics{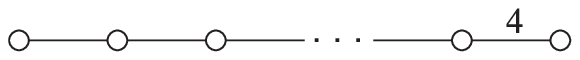}}
\end{array}$ \ \ ($l\geq 2$ vertices)

\medskip

$D_l$ \ \
$\begin{array}{c}
\mbox{\includegraphics{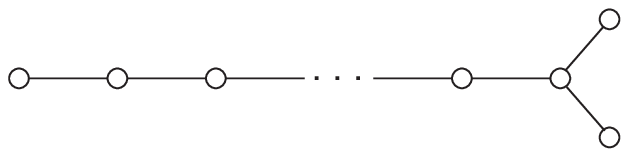}}
\end{array}$ \ \ ($l\geq 4$ vertices)

\medskip

$E_6$ \ \ $\begin{array}{c}
\mbox{\includegraphics{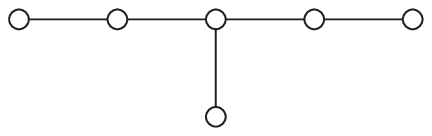}}
\end{array}$

\medskip

$E_7$ \ \
$\begin{array}{c}
\mbox{\includegraphics{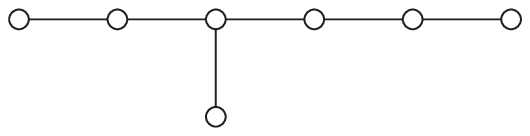}}
\end{array}$

\medskip

$E_8$ \ \
$\begin{array}{c}
\mbox{\includegraphics{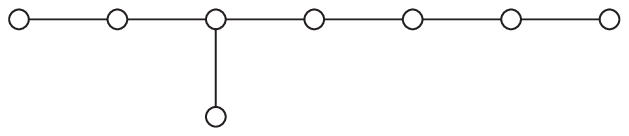}}
\end{array}$

\medskip

$F_4$ \ \ $\begin{array}{c}
\mbox{\includegraphics{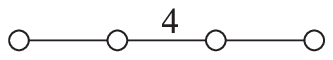}}
\end{array}$

\medskip

$G_2$ \ \ $\begin{array}{c}
\mbox{\includegraphics{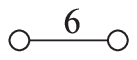}}
\end{array}$

\medskip

$H_3$ \ \ $\begin{array}{c}
\mbox{\includegraphics{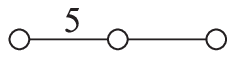}}
\end{array}$

\medskip

$H_4$ \ \ $\begin{array}{c}
\mbox{\includegraphics{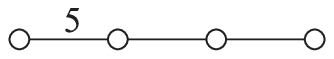}}
\end{array}$

\medskip

$I_2(p)$ \ \ $\begin{array}{c}
\mbox{\includegraphics{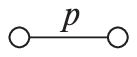}}
\end{array}$ \ (either $p=5$ or $p\geq 7$)

\medskip

\noindent and if $(W,S)$ is an irreducible Coxeter system with finite
set $S$, then the associated quadratic form is positive and
generated iff Coxeter graph is isomorphic to one of $f$-graphs of
the list IV:

\medskip

$\widetilde{A_1}$ \ \ $\begin{array}{c}
\mbox{\includegraphics{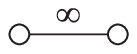}}
\end{array}$

\medskip

$\widetilde{A_l}$ \ $(l\geq 2)$ \ \
$\begin{array}{c}
\mbox{\includegraphics{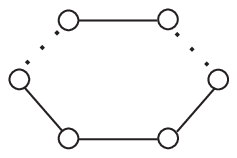}}
\end{array}$
\ \ (cycle with $l+1$ vertices)

\medskip

$\widetilde{B_2}$ \ \ $\begin{array}{c}
\mbox{\includegraphics{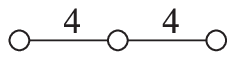}}
\end{array}$

\medskip

$\widetilde{B_l}$ \ $(l\geq 3)$ \ \
$\begin{array}{c}
\mbox{\includegraphics{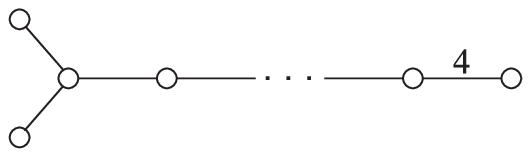}}
\end{array}$ \ \ ($l+1$ vertices)

\medskip

$\widetilde{C_l}$ \ $(l\geq 3)$ \ \
$\begin{array}{c}
\mbox{\includegraphics{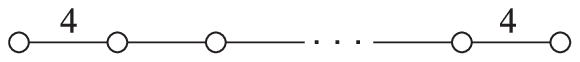}}
\end{array}$ \ \ ($l+1$ vertices)

\medskip

$\widetilde{D_l}$ \ $(l\geq 4)$ \ \
$\begin{array}{c}
\mbox{\includegraphics{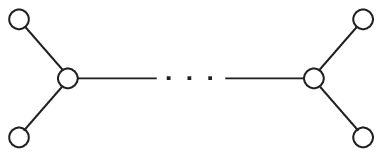}}
\end{array}$ \ \ ($l+1$ vertices)

\medskip

$\widetilde{E_6}$ \ \
$\begin{array}{c}
\mbox{\includegraphics{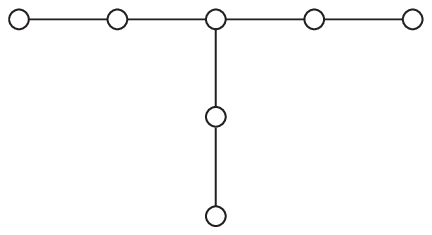}}
\end{array}$

\medskip

$\widetilde{E_7}$ \ \
$\begin{array}{c}
\mbox{\includegraphics{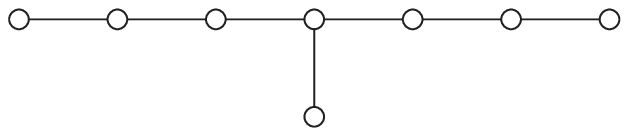}}
\end{array}$

\medskip

$\widetilde{E_8}$ \ \
$\begin{array}{c}
\mbox{\includegraphics{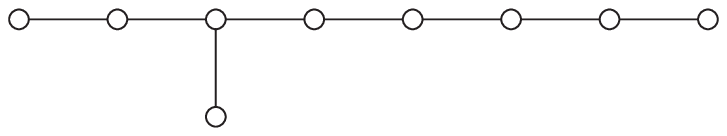}}
\end{array}$

\medskip

$\widetilde{F_4}$ \ \
$\begin{array}{c}
\mbox{\includegraphics{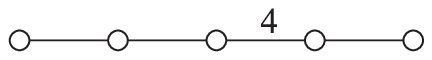}}
\end{array}$

\medskip

$\widetilde{G_2}$ \ \ $\begin{array}{c}
\mbox{\includegraphics{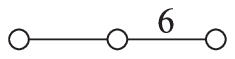}}
\end{array}$

\medskip

Recall that for Coxeter graphs $f(\alpha)\geq 3$, and
$f(\alpha)=3$ if no number is written over an edge.

For a Coxeter graph $\overline{\Gamma} =(\Gamma,f)$
we denote $\widehat{f}$-graph $(\Gamma,\widehat{f})$, where
$\widehat{f}(\alpha)=4\cos^2\dis{\pi\over f(\alpha)}$,
by $\widehat{\Gamma}$.
Of course, $\widehat{\Gamma}$ is not a Coxeter graph,
and, in general, it is not an integral $\widehat{f}$-graph.
For $f(\alpha)=\infty$ we assume $\widehat{f}(\alpha)=4\cos^2(0)=4$.
Denote the $\rho$-degree in $\hat \Gamma$ by $\hat g_\rho$
and the value $\displaystyle \frac{\partial v}{\partial \alpha}$ in
$\hat\Gamma$ by $\displaystyle \widehat{\frac{\partial v }{\partial \alpha}}$.

\begin{proposition}
A connected Coxeter graph $\overline{\Gamma}$ belongs to III
(resp. IV) iff for any point $x$ of the graph $\widehat{\Gamma}$
$\rho$-degree $\hat g_\rho(x)$ is less than 4
(resp. is less or equal to 4,
and there is $y$ such that $\hat g_\rho(y)=4$).
\end{proposition}

If $\widehat{\Gamma}$ is an integral $\widehat{f}$-graph
(i.e. $f(\alpha)\in\{3,4,6,\infty$ for $\alpha\in \Gamma_e\}$),
then the statement (in both sides) follows from proposition~2 and
a comparison III and \underline{I}, and IV and~\underline{II}.

Suppose that $\overline{\Gamma}$ is an $f$-graph satisfying
the conditions of the proposition such that $\hat \Gamma$ is not
integral. Note that if $\widehat f(\alpha)\not=1$, then $\widehat
f(\alpha)\geq 2$.

If $6<f(\alpha)<\infty$, then $\widehat{f}(\alpha)>3$;
from the definition of $\rho$-degree it follows that
$\varphi(\alpha)=\{x,y\}$, $g(x)=1$, $g(y)=1$,
and we have $I_2(p)$ ($p>6$).

Thus, it remains to consider only the case
\begin{gather*}
\alpha\in \Gamma_e,\quad f(\alpha)=5,\quad
\varphi(\alpha)=\{x,y\},\quad
\widehat{f}(\alpha)=4\cos^2\dis{\pi\over 5}=\\
=4\left({1+\sqrt{5}\over 4 }\right)^2={(1+\sqrt{5})^2\over
4}=\widehat\alpha>\mbox{2.5}.
\end{gather*}
$\hat g_\rho (x)\leq 4$ implies $g(x)\leq 2$.
If $g(x)=g(y)=1$, then $\overline{\Gamma}=I_2(5)$. Let
\begin{gather*}
g(x)=2,\quad \varphi(\beta)=(z,x), \quad \widehat
f(\beta)=\widehat\beta,
\\ \widehat g_\rho(x)=\widehat{\beta} \widehat{{\p x\over \p
\beta}}+\widehat{\alpha}\widehat{{\p x\over \p \alpha}}\leq 4,
\end{gather*}
 $\widehat{\beta}=1$ (if $\widehat{\beta}>1$, then$\widehat{\beta}\geq 2$,
 $\widehat{g}_\rho(x)>4$), $\dis \widehat{{\p x\over \p \alpha}}=1$;
$\dis \widehat{{\p x\over \p \beta}}<\mbox{1,5}=\rho(3)$, the pair
$(x,\beta)$ is $\widehat{f}$-simple, $|\Gamma'(x,\beta)|<3$.
Hence, only $H_3$ and $H_4$ are possible. On the other hand,
$\widehat{\alpha}<2 \dis{2\over 3}$, $\rho(2)=\dis 1\frac 13$,
therefore $H_3$, $H_4$ (and $I_2(p)$) satisfy the conditions of
the proposition.

\begin{remark}\rm
It is known that the list III characterizes finite groups generated by
reflections.

Let $G$ be a group generated by formatives $a_1,\ldots,a_t$
and relations $a_i^2=1$, $(a_i a_j)^{n_{ij}}=1$. A pair $a_i$,
$a_j$ is {\it special} if $n_{ij}>3$. Denote the number of
elements $a_j$ not commutating with $a_i$ by $g(a_i)$.
Subset $H\subset A=\{a_1,\ldots,a_t\}$ is
a {\it component} if $H$ cannot be represented in the form $H_1\cup H_2$
such that $a_i a_j=a_j a_i$ for all $a_i\in H_1$, $a_j\in H_2$.
Any $X\subset A$ has a unique representation as a union of disjoint components
$X_1,\ldots, X_{c(X)}$ such that if $a_i\in X_k$, $a_j\in X_l$ and
$k\not= l$, then $a_ia_j=a_ja_i$. We will assume $c(A)=1$
(that corresponds to connectivity of a Coxeter graph).
It immediately follows from the list III that a group  $G$ is finite iff
$G$ is subjected to one of the following conditions:

A) $A$ does not contain special pairs, contains $a_r$ such that
$g(a_r)=1$  and contains at most one generator $a_i$ such that $g(a_i)>2$.
In this case if $g(a_i)>2$, then
$g(a_i)=c(X)$, where $X=A\backslash a_i$ and $\rho(|X_1|,\ldots,|X_{g(a_i)}|)<4$;

B) $A$ contains precisely one special pair  $(a_i,a_j)$,
$g(a_s)\leq 2$ for $i=\overline{1,t}$  and there exists $a_r$
such that $g(a_r)=1$. Moreover,

1) if $n_{ij}\geq 6$, then $t=2$,

2) if $n_{ij}=5$, then $t\leq 4$ and $\min\{g(a_i),g(a_j)\}=1$,

3) if  $n_{ij}=4$, then either  $t\leq 4$, or $\min\{g(a_i),g(a_j)\}=1$.
\end{remark}

\section{Marked quivers}

In section~5 and further the matrix problems become the main object
of our study~\cite{11}. From naive point of view these are the
problems of an equivalence of matrices or sets of matrices by
means of some admissible transformations. For example,
representations of a poset $S$ are matrices divided into $n$
($=|S|$) vertical bands, one for each element of $S$,
and it is admitted to add columns of the band $s_i$ to columns
of $s_j$ iff $i\leq j$, and it is also admitted to make any elementary
transformations of rows (and columns into the same band).

On the other hand, such problems can be naturally formulated
in categorical terms.

In this section we remind some known facts and introduce a
categorical terminology applying further.

Note two points where our terminology and generally accepted one
are different.

Firstly, we use the right notation: if $\gamma:X\longrightarrow Y$,
$\delta:Y\longrightarrow Z$, then we denote the product of $\gamma$ and $\delta$
by $\gamma\delta$ (not by $\delta\gamma$).

Secondly, we consider modules (see~\cite{15}) and bimodules
over categories (generalizing them even for non-semiadditive categories)
although instead of them it is always possible to
consider functors and bifunctors. However,
in classical representation theory
(for instance, of finite-dimensional algebras)
it also would be possible to consider representations,
i.e. homomorphisms to matrix algebras, instead of modules,
but a consideration of modules is preferable.

{\it The category of morphisms} $K^\triangle$ of an arbitrary
category $K$ is called (\cite{16}) a category defined in the following way:
$\Ob K^\triangle=\Mor K$, and if $\varphi,\varphi'\in\Mor K$,
$\varphi:A\longrightarrow B$, $\varphi':A'\longrightarrow B'$
$(A,B,A',B'\in\Ob K)$, then
$K^\triangle(\varphi,\varphi')=\{(\alpha,\beta)\ |\ \alpha\in
K(A,A'), \beta\in K(B,B')/ \varphi\beta=\alpha\varphi'\}$, i.e.
the diagram below is commutative
\begin{equation}
\begin{CD}
A @>{\varphi}>> B\\
@V{\alpha}VV @VV{\beta}V\\
A' @>{\varphi'}>> B'
\end{CD}\tag{*}
\end{equation}

A natural generalization of this notion is {\it representations of
quiver}~\cite{10}.

Let $Q$  be a (finite) {\it quiver} (i.e. oriented graph). $Q_v$
is the set of its vertices and $Q_a$ of arrows of $Q$.
The beginning (tail) and the end (head) of an arrow $\alpha$
are denoted by $t_\alpha$ and by $h_\alpha$ respectively.
A {\it representation} $T$ of a quiver $Q$ in a category $K$
associates an object $T(x)\in\Ob K$ for each $x\in Q_v$,
and a morphism $T(\alpha):T(t_\alpha)\longrightarrow T(h_\alpha)$
for each $\alpha\in Q_a$.

A morphism from representation $T$ to representation $T'$ in
the {\it category $K^Q$ of representations of a quiver} $Q$ is a set of
morphisms $\alpha_x:T(x)\longrightarrow T'(x)$, one for each $x\in Q_v$,
such that for any $\alpha\in Q_a$ the following diagram is commutative
\begin{equation}
\begin{CD}
T(x) @>{T(\alpha)}>> T(y)\\
@V{\alpha_x}VV @VV{\alpha_y}V\\
T'(x) @>{T'(\alpha)}>> T'(y)
\end{CD}\tag{**}
\end{equation}

A product of morphisms is naturally defined for categories of morphisms
as well as for categories of representations of quivers.
Thus, a category of morphisms $K^\triangle$ can be considered
as a category of representations of quiver
$\triangle:a\overset{\alpha}\longrightarrow b$ ($T_v=\{a,b\}$,
$T_a=\{\alpha\}$, $h_\alpha=b$, $t_\alpha=a$).

\begin{remark}\rm Usually representations of quivers are considered in
the category $\mod k$, where $k$ is a field or a commutative ring.
We do not require additivity and even semiadditivity of
the category $K$. We give an example of non-additive category, where
representations of quivers seemingly have a substantial interest.
Let $H$ be {\it a category of Hilbert spaces}, where
$H(A,B)=\{\varphi\in\Hom_c(A,B)\ |\ \varphi\varphi^*=1_A\}$.
In fact, representations of quivers and quivers with relations in $H$
are considered in~\cite{25}, moreover, the construction of
Coxeter functors~\cite{17} extends to this case that yields
a progress in several problems of functional analysis.
\end{remark}

For each quiver $Q$, it can be naturally associated
a (non-oriented) graph $\Gamma(Q)$.
The criterion of finite representativity of a quiver was obtained in~\cite{10}
(see also~\cite{17}), and the tameness criterion (independently) in~\cite{19,18}
(see also~\cite{15}).
A connected quiver $Q$ is finitely represented (resp. tame) iff
$\Gamma(Q)$ is a scheme (resp. extended scheme) of Dynkin without
multiple edges, i.e.  $A_e$, $D_e$, $E_6$, $E_7$, $E_8$
(resp. $\widetilde{A_e}$, $\widetilde{D_e}$, $\widetilde{E_6}$,
$\widetilde{E_7}$, $\widetilde{E_8}$).

A replacement of the category $K$ by a bimodule (over two categories)
is another natural generalization of a category of morphisms.

We say that $M_K$ is a (right) {\it moduloid} over (possible not
semiadditive) category $K$, if a set  $M(A)$ is attached to
each $A\in\Ob K$, and a map from $M(A)$ to $M(B)$ is attached to
each $\varphi\in K(A,B)$, and for $a\in M(a)$ and $\varphi\in M(A,B)$
is defined a product $a\varphi\in M(B)$ such that
\begin{enumerate}
\item
$a1_A=a(a\in A\in\Ob K)$
\item
$a(\varphi\psi)=(a\varphi)\psi$, ($\varphi\in K(A,B)$, $\psi\in
K(B,C)$).
\end{enumerate}

If $K$ is $k$-category ($k$ is a field or a commutative ring)
\cite{15}, then $M_k$ is a {\it $k$-linear module} over $K$ if each
$M(A)$ is a finitely generated module over $k$, and besides 1 and 2,
\begin{enumerate}
\item[3.]
$(\lambda a+\mu b)(\alpha\varphi+\beta\psi)=
(\lambda\alpha)a\varphi+(\lambda\beta)a\psi+
(\mu\alpha)b\varphi+(\mu\beta)b\psi$,
\end{enumerate}
holds, where $\lambda,\mu,\alpha,\beta\in k$; $a,b\in M(A)$;
$\psi,\varphi\in K(A,B)$; $A,B\in\Ob K$.

Similarly, one can define a left moduloid and a left $k$-linear module ${}_KM$
as well as a bimoduloid and a $k$-linear bimodule ${}_KM_L$ over two,
in general, different categories $K$ and $L$ (in the last two cases
it is given sets and finitely generated $k$-module $M(A,B)$, where
$A\in\Ob K$ and $B\in\Ob L$).

We call the elements of sets $M(A)$ $(A\in\Ob K)$
for left (right) moduloids and $k$-linear modules and sets $M(A,B)$
for moduloids and $k$-linear bimodules
{\it representatives} of a moduloid, $k$-linear module,
bimoduloid and $k$-linear bimodule respectively.

A moduloid is {\it faithful} if $\{\varphi,\psi\}\subset K(A,B)$
implies $M(A)\ni a$, $a\varphi\not= a\psi$.

A {\it category of representatives} of $M^\triangle$ bimoduloid
$M={}_KM_L$ is a category whose objects are representatives of
$M$, and if $\varphi\in M(A,B)$, $\varphi'\in M(A',B')$, then
$M^\triangle(\varphi,\varphi')=\{(\alpha,\beta)\ |\ \alpha\in
K(A,A'),\beta\in L(B,B'),\varphi\beta=\alpha\varphi'\}$, i.e. the
diagram ($*$) is commutative. Thus, this definition literally
repeats the definition of a category of morphisms.

If $K'$ is a subcategory of an arbitrary category $K$, then
it is defined the bimoduloid ${}_KM(K')_{K'}$ (over $K$ and $K'$),
where $M(K')(A,B)=\Hom_K(A,B)$ $(A\in\Ob K$, $B\in\Ob K')$.

If $K'$, $K''$ are two subcategories of a category $K$, then
analogously it is defined the bimoduloid ${}_{K'}M(K',K'')_{K''}$,
where $M(A,B)=\Hom_K(A,B)$ $(A\in\Ob K'; B\in\Ob K'')$.

Further we assume that $K=\mod k$, $k$ is an algebraically
closed field, $K'$ is a subaggregate of $K$ (i.e. an additive
$k$-subcategory whose idempotents are splitted~\cite{15}). Under these
assumptions we denote the category $M(K')^\triangle$ by $\Rep K'$.
An arbitrary faithful $k$-linear module over a $k$-aggregate can be
identified as a subaggregate of $K$. Further saying ``aggregate''
we mean a subaggregate of~$K$.

The category $\Rep K'$ was introduced in \cite{20} and considered
in \cite{21,15}. This category plays an important role in the theory
of representations. In particular, in \cite{22} it is shown that
the category of representations of a finite dimensional algebra
can be ``reduced'' to it.

Combining two given generalizations of a category of morphisms
we come to the notion of a representation of marked quiver.

We say that a quiver $Q$ is {\it marked} if to each $a\in Q_v$
is attached a category $K(a)$, and to each $\alpha\in Q_a$ is
attached a bimoduloid ${}_{K(a)}(M^\alpha)_{K(b)}$, where
$t_\alpha=a$, $h_\alpha=b$ (possibly, $a=b$).

A representation $T$ of a marked quiver $\overline{Q}$ attaches
an object $T(a)\in K(a)$ to each $a\in Q_v$ and
a representative $T(\alpha)\in M^\alpha(T(x),T(y))$ to each $\alpha\in Q_a$,
where $x=t_\alpha$, $y=h_\alpha$.

A morphisms from a representation $T$ to a representation $T'$
of {\it category $\Rep(\overline{Q})$ of representations of a marked quiver
$\overline{Q}$} is a set of morphisms $\alpha_x\in K(x)(T(x), T'(x))$,
one for each $x\in Q_v$, such that for any $\alpha\in Q_a$ the diagram ($**$)
(in the bimoduloid ${}_{K(x)}(M^\alpha)_{K(y)}$!) is commutative
($x=t_\alpha$, $y=h_\alpha$).

A marked quiver $\overline{Q}$ (resp. bimodule $M$) is {\it finitely
represented} if the category $\Rep{\overline{Q}}$ (resp.
$M^\triangle$, which can be considered as a partial case
of $\Rep \overline{Q}$ assuming $Q=\triangle$) has a finite
number of indecomposable isoclasses.

We assume that the quiver $Q$ is connected and
$Q_a\not=\varnothing$. We also assume that the marked quiver
$\overline{Q}$ is {\it $k$-marked}, i.e. all $K(a)$ $(a\in Q_v)$
are subaggregats of $K$, and each ${}_{K(x)}M_{K(y)}^\alpha$
($\alpha\in Q_\alpha$, $t_\alpha=x$, $h_\alpha=y$) is
$M(K(x),K(y))$. So in this case bimodules $M^\alpha$ $(\alpha\in
Q_a)$ are determined uniquely by aggregates $K(a)$ $(a\in Q_v)$.
Representations of marked, but not $k$-marked quivers, and
representations of $k$-marked {\it quivers with relations} are
also interesting, but will not be considered in this article.

Most of matrix problems can be considered as
representations of  $k$-marked quivers. However, to conceive
categories $\Rep{\overline{Q}}$, it is necessary to choose an
appropriate language for consideration of subaggregats of the category $K=\mod k$.
Among these subaggregats, ones generated by posets as well as posets with
an equivalence relation and biequivalence relation play an important role.

A poset with an equivalence $\overset{\sim}{S}$ is a (finite) set
$S$ with two (completely independent) relations given on $S$:
a partial order $\leq$ and an equivalence~$\sim$.

A {\it biequivalence} given on a (finite) poset $S$
is an equivalence relation $\approx$ on $\{(s,t)\ |\ s,t\in S,s\leq t\}$
such that
\begin{enumerate}
\item[i)]
if $(s_1,t_1)\approx(s_2,t_2)$ and $(s_1,t_1)\ne(s_2,t_2)$, then
$s_1\ne s_2$, $t_1\ne t_2$;
\item[ii)]
if $(s_1,t_1)\approx(s_2,t_2)$ and $s_1\leq x_1\leq t_1$, then
there exists $x_2$ such that $s_2\leq x_2\leq t_2$,
$(s_1,x_1)\approx(s_2,x_2)$ and $(x_1,t_1)\approx(x_2,t_2)$.
\end{enumerate}
From these conditions it follows that $x_2$ is uniquely defined,
and the following conditions hold:

\begin{enumerate}
\item[iii)] if $(a,a)\approx (c,d)$, then  $c=d$ (setting $a=x_1$,
obtain from ii) $(a,a)\approx (c,x_2)\approx (c,d)$, and from i) $c=d$);

\item[iv)] if $(s_1,t_1)\approx (s_2,t_2)$, then
$(s_1,s_1)\approx (s_2,s_2)$, $(t_1,t_1)\approx (t_2,t_2)$
(setting $x_1=s_1$ we get $(s_1,s_1)\approx (s_2,x_2)$, and
from iii) $x_2=s_2$).
\end{enumerate}

We denote a poset $S$ with a biequivalence by $\overset{\approx}{S}$.
The relation $\approx$ induces an equivalence relation $\sim$
on $S$ itself ($s\sim t$ if $(s,s)\approx(t,t)$). On the other hand,
$\overset{\sim}{S}$ can be considered as a partial case of $\overset{\approx}{S}$
letting for $a,b,c,d\in S$, $ (a,b)\approx(c,d)$ if (and only if) $a= b$,
$c= d$ and $a\sim c$.

A biequivalence is {\it transitive} if $(s_1,t_1)\approx(s_2,t_2)$
and  $(t_1,u_1)\approx(t_2,u_2)$ imply
$(s_1,u_1)\approx(s_2,u_2)$. We write $s\nc t$ if $s$ and $t$
($\in(S,\leq)$) are not comparable.

For $\overset{\approx}{S}$, one can associate
a subaggregate $K(\overset{\approx}{S})$ in $K$.

First, construct a full subcategory $L_X$ of {\it
$X$-graduated} spaces of $K$ for an arbitrary finite set $X$.

An object $V\in L_X$ is $\bigoplus\limits_{x\in X}V_x$, where $V_x\in
K$. Naturally, if $\psi\in L_X(V,W)$, then $\psi=\sum\limits_{x,y\in
X}\psi_{x,y}$, where $\psi_{x,y}\in K(V_x,W_y)$.

Aggregate $K(\overset{\approx}{S})$ is a subcategory of $L_S$
given by the following conditions
\begin{enumerate}
\item
$V\in\Ob K(\overset{\approx}{S})$, if $x\approx y$ implies
$V_x=V_y$ $(x,y\in S)$;
\item
$\psi\in K(\overset{\approx}{S})(V,W)$, if $x\nc y$ or $x>y$
implies $\psi_{x,y}=0$, and $(x,y)\approx(u,v)$ implies
$\psi_{x,y}=\psi_{u,v}$ ($x,y,u,v\in S$).
\end{enumerate}

In particular, if $\approx$ is trivial or not trivial only on pairs
$(s,s)$ (i.e. actually we have an equivalence $\sim$ on the poset $S$),
we obtain aggregates $K(S,\leq)$ and $K(\overset{\sim}{S})$.
The category $K(S,\leq)^\triangle$ coincides with $\Rep(S,\leq)$,
introduced in \cite{3}, $K(\overset{\sim}{S})$ coincides with the
category of representations of posets with an
equivalence~\cite{23}. We also denote the categories
$K(\overset{\approx}{S})^\Delta$ and $K(\overset{\sim}{S})^\Delta$
by $\Rep \overset{\approx}{S}$ and $\Rep \overset{\sim}{S}$, and
the sets of indecomposable isoclasses of these categories by
$\ind\overset{\approx}{S}$ and $\ind\overset{\sim}{S}$ respectively.
An aggregate $K(\overset{\approx}{S})$ is {\it transitive}
if biequivalence $\approx$ is transitive.

From matrix point of view the category $\Rep\overset{\approx}{S}$
(strongly speaking, up to an equivalence)
can be given in the following way. Every its
objects is a matrix $T$ divided into $n$ ($=|S|$) vertical bands
$T_s$ ($s\in S$). $\Hom(X,Y)$ consists of pairs of matrices
$(A,B)$ such that $AX=YB$, and $B$ is divided into blocks $B_{st}$
($s,t\in S$) according to the divisions of matrices $X$ and $Y$,
and, moreover, $B_{st}\not=0$ only if $s\leq t$ and $B_{st}=B_{uv}$,
when $(s,t)\approx (u,v)$.

Note that in some cases the categorical language is more preferable,
whereas in other case the matrix language is more suitable.
For example, the following statement is obvious in matrix language.

\begin{lemma}
Let $S=U\cup V$, $U\cap V=\varnothing$, Suppose that $(ab)\approx (cd)$
implies either $a,b,c,d\in U$, or $a,b,c,d\in V$, and $x\nc y$
implies either $x,y\in U$, or $x,y\in V$. Then $\ind
\overset{\approx}{S}=\ind\overset{\approx}{U}$ $\amalg\,
\ind\overset{\approx}{V}$.
\end{lemma}

A poset with an equivalence (or biequivalence) is {\it chain}
(rest. {\it antichain}) if $s\sim t$ implies $s$ and $t$
are comparable (resp. incomparable).

We denote the order of an equivalence class of $\sim$ containing $s$
by $\dim s$ $(s\in S)$. The dimension of $\overset{\approx}{S}$ is
$\max\dim s$, $ s\in S$. We call $\max\dim A$, where $A$ is
an indecomposable object $K'$, the dimension $\dim K'$ of
subaggregate $K'\subset K$. If $K'=K(\overset{\approx}{S})$,
then it is clear that $\dim K'=\dim \overset{\approx}{S}$.
The finite representativity and tameness for aggregates,
posets with an equivalence and biequivalence relations and
$k$-marked quivers are defined in the natural way.

Representations of $k$-marked quivers contain itself
representations of posets and representations of aggregates $\Rep
K'$) and also, of course, representations of (unmarked) quivers.
For posets and quivers, the criteria of finite representativity
and tameness are known (see sections~2,~5). Obtaining of analogous
criteria for $k$-marked quivers is apparently one of the main
problems of the theory of matrix problems. The main difficulties
here are contained in the partial case of representations of
aggregates.

The problem of its finite representativity is considered in
\cite{4,5,26}, see also \cite{15,27,28}. It is well known that
if an aggregate $\mathcal A$ is finitely represented, then
$\dim {\mathcal A}\leq 3$, and if ${\mathcal A}=K(\overset{\approx}{S})$,
then $\overset{\approx}{S}$ is chain, whereas if ${\mathcal A}$ is
tame, then $\dim {\mathcal A}\leq 4$. For dimension 2 the criterion
of finite representativity is given in \cite{4}, and for dimension 3
such criterion is formulated in~\cite{5}, where, however,
only the necessity of the given conditions is proved.

A representation of an aggregate of dimension 1 is a
representation of a poset. Representations of
$K(\overset{\approx}{S})$ of dimension 2 are considered in
section~7, and  representations of $K(\overset{\sim}{S})$ in
section~6.

For  $\overset{\approx}{S}$ set  $S^2=\{(x,y)\;|\; x,y\in S,
x<y\}$. From iii  it follows that the relation $\approx$
disintegrates into a relation $\sim$ on $S$ and an equivalence
relation on $S^2$. If $(x,y)\in S^2$, then we call the order of
the equivalence class $\approx$ in $S^2$ containing $(x,y)$ the
{\it the rank} ${\rm rank}\, (x,y)$ of the pair $(x,y)$, and let
${\rm rank}\,\overset{\approx}{S}=\max\,{\rm rank}\,(x,y)$. If
(and only if) ${\rm rank}\,\overset{\approx}{S}=1$, then $\Rep
\overset{\approx}{S}\simeq \Rep \overset{\sim}{S}$. Clearly, ${\rm
rank}\, \overset{\approx}{S}\leq \dim \overset{\approx}{S}$.

We write $x\lhd x$, if $\dim x=1$, $x\lhd y$, if ${\rm
rank}\,(x,y)=1$,
 and $x\Rightarrow y$, if ${\rm rank}\,(x,y)> 1$.

From i and ii it follows that if either $a\leq b\lhd c$, or $a\lhd
b\leq c$, then $a\lhd c$.

{\it Biordered set} (boset) $S^\lhd$ is a poset $S$ with
an additional relation $\lhd$ such that

1) if $a\lhd b$, then $a\leq b$;

2) if either $a\lhd b\leq c$, or $a\leq b\lhd c$, then $a\lhd c$.

If $\dim\overset{\approx}{S}=2$, then the relation $\approx$ is
uniquely determined by the relations $\lhd$ and $\sim$.
We call the elements of $S$ {\it points}.
A point $s$ is {\it small} if $\dim s=1$ and
{\it big} if $\dim s>1$. We denote small points by $\circ$,
and big points by $\bullet$; $\overset{\circ}{S}$ (resp.
$\overset{\bullet}{S}$) is the set of small (resp. big) points of $S$.
If $\dim s=2$, then by $s^*$ we denote such (unique) point that
$s^*\sim s$ and $s^*\ne s$. If  $U,V\subset S$,
then $U\nc V$ means that $u\nc v$ for any  $u\in U$,
$v\in V$. If $X,Y\subset S$, then $X^\nc(Y)=\{x\in X\ |\ \{x\}\nc
Y\}$ (if $y\in S$, then $X^\nc(y)=X^\nc(\{y\})$).

A 1-{\it chain} of a poset with equivalence relation
$\overset{\approx}{S}$ (or boset $S^\lhd$) is
$Z=\{s_1,\ldots,s_t\}\subset S$ $(t\geq 0)$, if $s_i\lhd s_j$ for
$1\leq i<j\leq t$. (For $\overset{\sim}{S}$ any chain is a 1-chain).

Introduce the notation of {\it normality} for big points.
A point $t$ is {\it 1-normal}, if $S^\nc(t)$ is a chain of $\overset{\circ}{S}$.

Further, if $j$-normality is already defined for $j<i$, then a point
$t$ is {\it $i$-normal}, if
\begin{enumerate}
\item
$S^\nc(t)$ is a 1-chain;
\item
if $x\in S^\nc(t)$, then $\dim x\leq 2$, and if $\dim x\ne 1$,
then $x^*$ is $j$-normal for some $j<i$.
\end{enumerate}

A point is {\it normal} if it is $i$-normal for some $i$. A point $y$
is {\it conormal} if $\dim y=2$ and $y^*$ is normal; here if $y^*$
is $i$-normal, then $y$ is $i$-{\it conormal}. Note that some
big points can be normal and conormal, whereas others are neither normal,
nor conormal. Points of dimension more than 2 can be normal but
cannot be conormal.

\begin{example}\rm
$\overset{\approx}{S}=\overset{\sim}{S}=$ \ $\begin{array}{c}
\mbox{\includegraphics{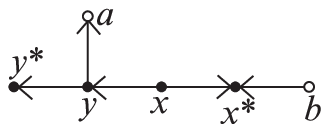}}
\end{array}$
Here $x$ is 1-normal point $(S^\nc(x)=\{b\})$, $y$ is 2-normal
point $(S^\nc(y)=\{b,x^*\})$, $x^*$, $y^*$ are conormal but not
normal $(S^\nc(x^*)=\{y,y^*,a\}$, $S^\nc(y^*)=\{a,x^*,b\})$.
\end{example}

A poset $S$ is a {\it $p$-poset} $(S,p)$,
if it is given a function $p$ on $S$ with the values in ${\mathbb N}^\infty$.
$(\overset{\approx}{S},p)$ denotes $p$-poset with a biequivalence
relation. In this case, for $Z\subset S$ set
$p(Z)=\sum\limits_{z\in Z} p(z)$ if $Z$ is 1-chain,
and $p(Z)=\infty$, if $Z$ is not 1-chain.

Construct on $\overset{\approx}{S}$ a function
$\overset{\approx}{p}$ uniquely defined by the relation $\approx$.

First, note that if $x\in \overset{\approx}{S}$ is an $i$-conormal
point, and $S^\nc(x^*)=Y$, then  $Y$ is 1-chain, and if $y\in Y$,
then $y$ is small or a $j$-conormal point, where $j<i$. If $\dim
s=1$, then we set $\overset{\approx}{p}(s)=1$. If $\dim s>1$ and
$s$ is not conormal, then we set $\overset{\approx}{p}(s)=\infty$.
Further, we define $\overset{\approx}{p}$ for conormal points
consecutively by the degree of its conormality. If the function
$\overset{\approx}{p}$ is defined for all $j$-conormal points for $j<i$,
and $x$ is $i$-conormal, then we set
$\overset{\approx}{p}(x)=2+\sum\limits_{y\in
Y}\overset{\approx}{p}(y)$ ($Y=S^\nc(x^*)$).

For $\overset{\sim}{S}$, one can construct a function $\overset{\sim}{p}$
using usual chains instead of 1-chains in its definition.

In the example 3 $\overset{\sim}{p}(x^*)=3$,
$\overset{\sim}{p}(y^*)=6$,
$\overset{\sim}{p}(x)=\overset{\sim}{p}(y)=\infty$.

\section{Representations of posets with equivalence}

Representations of posets with equivalence were considered in
\cite{1,29,23} (in \cite{1} they were called ``representation of
weakly completed posets''). In  \cite{1} the problems of finite
representativity of a poset with equivalence $\overset{\sim}{S}$
and the description of its representations were obtained by
reduction to representations of posets.

In \cite{29} a tameness criterion was given for antichain (if $x\sim y$,
then $x\nc y$) posets with equivalence. Namely, it is proved that infinitely
represented antichain poset $S$ is tame iff it does not contain the following
critical subsets $N_0-N_5$ (see section~2) in $\overset{\circ}{S}$ and

$N_6$ \ \ $\bullet$ \ \ $\bullet$ \ \ $\circ$;

$N_7$ \ \ $\bullet$ \ \ $\circ$ \ \
\raisebox{2.5mm}{$\begin{array}{c}\circ\\[-2.5mm] |\, \\[-2.5mm] \circ\end{array}$};

$N_8$ \ \ $\bullet$ \ \ $\bullet$ \ \ $\bullet$;

$N_9$ \ \ $\bullet$ \ \ $\circ$ \ \ $\circ$ \ \ $\circ$.

In \cite{23} it is proved that $N_0-N_9$\footnote{The list
$N_0-N_9$ was also announced by L.A. Nazarova as a tameness criterion
for antichain posets with a biequivalence relation and arbitrary
antichain  aggregates (in the last case it is given an exact
formulation what does ``inclusion'' of  $N_i$ in an aggregate mean),
see~\cite{30}.} are critical for a wider class of
``quasiantichains'' posets with equivalence (see below).
In general, for an arbitrary $\overset{\sim}{S}$,
there is no a direct tameness criterion in~\cite{23},
but it is shown how this question reduces to the given partial case.

In this section we give (in terms of $\rho$) direct criteria of
tame and finite representativity for an arbitrary
$\overset{\sim}{S}$, and in section~8 it is done for quivers
$k$-marked by subaggregates of the form $K(\overset{\sim}{S})$.

Note that the function $\rho$ defined on ${\mathbb N}$ will
correspond to $N_0-N_5$, whereas $\rho$ defined on ${\mathbb N}^\infty$
will correspond to $N_6-N_9$.

A set $\overset{\sim}{S}$ is {\it quasiantichain} if
$\dim \overset{\sim}{S}=2$ and $\overset{\sim}{S}$ does not contain
normal points (equivalently, $\overset{\sim}{S}$ does not contain 1-normal points,
i.e. each big point is included in  $\bullet\
\bullet$ or $\bullet \ \circ \ \circ$).
From \cite{1} it follows that a quasiantichain set cannot be
finitely represented.

Before to formulate finiteness and tameness criteria for $\overset{\sim}{S}$,
we go back to analogous questions for $(S,\leq)$ considered in section~2
and reformulate the criteria given there in the form suitable for generalization.

The corresponding criteria are formulated in section~2 in the following two forms:
\begin{enumerate}
\item
In form of absence of critical subsets $K_i$, $N_i$.
\item
In the form $\Rho(S)<4$, $\Rho(S)=4$, where $\Rho$ is a function defined on posets,
which gave birth to the function $\rho$.
\end{enumerate}
If $S$ is a primitive poset  $\bigsqcup\limits_{i=1}^tL^i$
($w(L^i)=1$), then $\Rho(S)=\rho(|L^1|\dots|L^t|)$. Set
$\rho_1(S)=\max\limits_{S'\subseteq S}\Rho(S')$, where $S'$ is a primitive poset.

It is clear that $\rho_1(S)\leq \Rho(S)$, but, for example, for
$S=\widehat{N}=\langle
2,2\rangle=$\raisebox{-1.5ex}{$\overset{\overset{\displaystyle
\circ }|}\circ$}$\displaystyle\diagdown$\raisebox{-1.5ex}
{$\overset{\overset{\displaystyle \circ }|} \circ $},
$\rho_1(S)=2{1\over 3}<\Rho(S)=2,4$. Call a poset $S$  {\it
quasiprimitive} if $S=\widehat{N}\sqcup Z$, $w(Z)\leq 1$.
$\Rho(\widehat{N}\sqcup Z)=2,4+\rho(|Z|)$.

$\Rho(\widehat{N}\sqcup Z)<4$ (resp. $=4$) iff
$|Z|<4$ (resp. $=4$), since $\rho(4)=1,6$.

For quasiprimitive $S=\widehat{N}\sqcup Z$ we let
$\rho_2(S)=\max|Z|$. For arbitrary $S$ we let
$\rho_2(S)=\max\limits_{S''\subset S} \rho_2(S'')$,
where $S''$ is quasiprimitive.
Finally, let $\rho(S)=\max\{\rho_1(S),\rho_2(S)\}$.

From the criteria in section~2 and proposition 1 it follows that
$S$ finitely represented (resp. tame) iff $\rho(S)<4$ (resp. $\rho(S)=4$).

Note that $\rho(S)\leq \Rho(S)$; here exact inequality takes place, in particular,
if $S$ is a uniform wattle (section~2) different from $\widehat{N}$
(for example, $\langle 3,3\rangle)$.

The criteria in terms of $\Rho$ are more natural, but in terms
of $\rho$ (then, actually, these criteria are equivalent to the absence of $K_i$, $N_i$)
they are more suitable for application.

Let $(S,p)$ be a $p$-poset. For a primitive
$S=\bigsqcup\limits_{i=1}^tL^i$ we set
$\rho'(S,p)=\sum\limits_{i=1}^t\rho(p(L^i))$; for
a quasiprimitive $S=\widehat{N}\sqcup Z$ we set
$\rho'(S,p)=\sum\limits_{z_i\in Z}p(z_i)$.
For an arbitrary $S$,
set $\rho(S,p)=\max\limits_{S'\subseteq S}\rho'(S',p)$,
where the maximum is taken over all primitive and quasiprimitive subsets
of $S'$ (and $S'=S$).

If $X\subset S$, then, of course, $\rho(X,p)$ is also defined.
Set $S_p=\{s\in S\;|\; p(s)<\infty\}$.

\begin{remark}\rm
If $\rho(S,p)<4$ and $A$ is an antichain of $S$, then
$|A|+|A^\infty|<4$, where $A^\infty=\{a\in A\; |\; a\not \in S_p\}$.
If for any antichain $A\subset S$, $|A|+|A^\infty|<4$ and
$\rho(S_p,p)<4$, then $\rho(S,p)<4$.

The first statement immediately follows from our definitions,
and the second from the fact that
if $S'$ is a primitive or quasiprimitive subset of $S$,
$\rho'(S',p)\geq 4$ and $S'\ni a\not\in S_p$,
then $a\in A^\infty$ and $|A|+|A^\infty|\geq 4$ for some antichain $A$.
\end{remark}

Set $\rho(\overset{\sim}{S})=\rho(S,\overset{\sim}{p})$.

Let $W=\{d_1,\ldots,d_t\}$ be a class of equivalence on $S$, $t>2$.
Set $\mu(W)=\sum\limits_{i=1}^t\rho(\overset{\sim}{p}(S^\nc(d_i))+1)$
and $\mu(\overset{\sim}{S})=\max\limits_{W\subset S} \mu(W)$.

\begin{proposition}
 $\overset{\sim}{S}$ is finitely represented iff
\begin{enumerate}
\item
$\rho(\overset{\sim}{S})<4$;
\item
$\mu(\overset{\sim}{S})<4$.
\end{enumerate}

$\overset{\sim}{S}$ is tame iff
$\max\{\rho(\overset{\sim}{S}), \mu(\overset{\sim}{S})\}=4$.
\end{proposition}

The proposition follows from the results  \cite{1,23},
we give a scheme of the proof.

First, note that for quasiantichain sets the tameness criterion
of proposition~6 coincides with the criterion in the form of absence $N_0$--$N_9$
proved in~\cite{23}.

If $\overset{\sim}{S}$ is quasiantichain, then $\overset{\sim}{p}(x)=1$
if $x$ is a small point, and $\overset{\sim}{p}(y)=\infty$ if $y$ is a big one.
It is easy to see that $\rho(\overset{\sim}{S})>4$, if $S\supset N_i$, $0\leq i\leq 9$.
On the other hand, it is easy to see that if
$S'$ is minimal (primitive or quasiprimitive) subset of $S$
such that  $\rho'(S',\overset{\sim}{p})>4$, then $S'$ has the form
$N_i$ $(0\leq i\leq 9)$.

If a set is not quasiantichain, then the point of our constructions
is the following: if $\overset{\sim}{p}(x)=t<\infty$,
in some sense we can replace $x$ by a chain consisting of $t$ small point
(for this, we throw out some ``not essential'' points).

If $x$ is a 1-normal point, then construct
$\overset{\sim}{S}{}'_x$ as follows: exclude the point $x$ from
$S$, and replace the point $x^*$ by
$t=\overset{\sim}{p}(x^*)=2+\overset{\sim}{p}(S^\nc(x))$ small
points $x_1^*<' x_2^*<'\cdots<' x_t^*$. On ``old'' points the
relations $\leq$ and $\sim$ are preserved, $x_i^*<'$ (resp. $>'$)
$y\in S$, if $x^*<$ (resp. $>$) $y$.

In example 3 $\overset{\sim}{S}{}'_x=\begin{array}{c}{}\\[-1mm]
 \mbox{\includegraphics{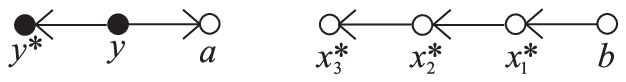}}
\end{array}$
The point $y$ became 1-normal.

The following lemma \cite{1}
(which can be easily proved in matrix language)
plays a crucial role in \cite{1} and
an important role in \cite{23}.

\medskip

\noindent {\bf Lemma.} {\it There is a natural one-to-one
correspondence between $\ind \overset{\sim}{S}$ and $\ind
\overset{\sim}{S}{}'_x$.}

\medskip

By $\overset{\sim}{p}{}'$ we denote the function constructed on
$\overset{\sim}{S}{}'_x$ analogically to $\overset{\sim}{p}$ for
$\overset{\sim}{S}$.

\begin{lemma}
If
$\max\{\rho(S,\overset{\sim}{p}),\rho(S'_x,\overset{\sim}{p}{}')\}\geq
4$, then
$\rho(S,\overset{\sim}{p})=\rho(S'_x,\overset{\sim}{p}{}')$.
\end{lemma}

$\overset{\sim}{p}(s)=\overset{\sim}{p}{}'(s)$ for any point $s \in S\cap S'_x$.
A point $x$ (1-normal) does not belong
to a quasiprimitive set or to a primitive $S'\subset S$
such that $\overset{\sim}{p} (S')\geq 4$;
$\overset{\sim}{p}(x^*)=\overset{\sim}{p}{}'\{x^*_1,\ldots,x^*_t\}=t$.

A poset with equivalence is called a {\it perfectly chain} if
\begin{enumerate}
\item[1)] $\dim  x=2$ implies that at least one of the points
$\{x,x^*\}$ is normal;

\item[2)] if $\dim y=3$, then at least two points of the equivalence class
$\{y,y',y''\}$ are 1-normal.
\end{enumerate}

In \cite{1} it is proved that finitely represented
$\overset{\sim}{S}$ is perfectly chain, and the same follows from
the conditions 1, 2 of proposition~5. Therefore, for $\dim S=2$
successively applying lemma~\cite{1}, it is possible to reduce
a finitely represented $\overset{\sim}{S}$ to a poset,
and the statement on finite representativity $\overset{\sim}{S}$ follows from
lemmas~\cite{1},~8 and the criterion of finite representativity of
posets reformulated above.

As it was noted, the tameness criterion for quasiantichain posets
with an equivalence relation coincides with the condition
$\rho(\overset{\sim}{S})= 4$. Hence, our tameness condition
for $\dim S=2$ follows from the result \cite{23} and lemmas
\cite{1} and 8, since applying these lemmas several times
an arbitrary $\overset{\sim}{S}$ can reduced to a quasiantichain poset.

If $\dim y=3$, $y\simeq u\simeq v$, and the points $u$, $v$ are 1-normal,
construct $(S''_y,\leq'',\sim'')$ excluding the points
$u$ and $v$ from $S$ and replacing $y$ by the poset $Y$,
which consists of small points $Y=U\times V$
(Descartes or {\it cardinal}~\cite{12} product),
$w(U)=w(V)=1$, $|U|=2+|S^\nc(u)|$, $|V|=2+|S^\nc(v)|$. If
$\overline{y}\in Y$, and $z\in S''_y\cap S$, then $\overline{y}<''
z$ (resp. $\overline{y} \ \; {}''\!\!> z$), if $y<z$ (resp.
$y>z$).

In this situation, similar to lemmas~\cite{1} and~8 statements
are proved in~\cite{1}.

Using these statements an arbitrary finitely represented
$\overset{\sim}{S}$ (taking into account that $\dim S\leq 3$ and $S$ is perfectly chain)
can be reduced to a poset.

For $\dim \overset{\sim}{S}=3$ the tameness criterion
is proved analogously. Note only that the condition 2) from the
definition of a perfectly chain poset with an equivalence relation
in the tame case, generally saying, does not take place as it is claimed
in~\cite{23}. However, this conditions is true if
there is no normal points of dimension 2
in~$\overset{\sim}{S}$. Hence, for $\dim \overset{\sim}{S}=3$,
first it is necessary ``to reduce'' all normal points of dimension 2
(``reducing'' $\overset{\sim}{S}$ to $\overset{\sim}{S}{}'_x$),
and then to reduce points of dimensional 3 (reducing
$\overset{\sim}{S}$ to $\overset{\sim}{S}{}''_y$), in the result
we obtain a quasiantichain set.

Similarly to lemma 8 one can prove that ``reducing''
if dimension of $x$ is 2 or dimension of $y$ is 3,
then $(S'_x,\leq',\sim')$ or $(S''_y,\leq'',\sim'')$
satisfies the condition $\overline{\mu}(\overset{\sim}{S}{}'_x)<4$ or
$\overline{\mu}(\overset{\sim}{S}{}''_y)<4$ (resp.
$\overline{\mu}(\overset{\sim}{S}{}'_x)=4$ or
 $\overline{\mu}(\overset{\sim}{S}{}''_y)=4$)
iff $\overline{\mu}(\overset{\sim}{S}_x)<4$ (resp.
$\overline{\mu}(\overset{\sim}{S}_y)=4$). In order to prove that
$\rho(\overset{\sim}{S}{}''_y)<4$ (resp.
$\rho(\overset{\sim}{S}{}''_y)=4$) iff
$\rho(\overset{\sim}{S})=4$, we use also the following
combinatorial statement, which can be easily checked (taking into
account the finiteness and tameness criteria, see~section~2).

\begin{lemma} If poset $S=(u-1)\amalg ((a+1)\times(b+1))$ $(u,a,b \in {\mathbb N})$,
then $S$ is finitely represented (resp. tame)  iff  $\rho(u,a,b)<4)$
(resp. $\rho(u,a,b)=4$).
\end{lemma}

If $\dim s=4$, then from tameness or the condition
$\mu(\overset{\sim}{S})\leq 4$ it follows that $S^\nc
(s)=\varnothing$, and this case can be reduced to the above one,
for instance, by lemma~7 (section~5).

\section{Finitely represented dyadic sets}

For $\overset{\approx}{S}$, two functions $\overset{\approx}{p}$ and
$\overset{\sim}{p}$ are defined. From their definitions
it follows that if $\overset{\approx}{p}(x)<\infty$, then
$\overset{\sim}{p}(x)=\overset{\approx}{p}(x)$, but if
$\overset{\approx}{p}(y)=\infty$, then, possibly,
$\overset{\sim}{p}(y)<\infty$.
Analogously to $\rho(\overset{\sim}{S})$ and $\mu(\overset{\sim}{S})$,
we define $\rho(\overset{\approx}{S})$ and
$\mu(\overset{\approx}{S})$ (replacing $\overset{\sim}{p}$ by
$\overset{\approx}{p}$).

\begin{proposition}
If $\overset{\approx}{S}$ is finitely represented, then
$\rho(\overset{\approx}{S})<4$ and
$\overline{\mu}(\overset{\approx}{S})<4$.
\end{proposition}

If $\overset{\approx}{S}$ is finitely represented or
$\mu(\overset{\approx}{S})<4$, then $\dim \overset{\approx}{S}\leq
3$. The statement follows from~\cite{7} if $\dim =1$,
from~\cite{4} if $\dim =2$ and from~\cite{5} if $\dim=3$
(if $\dim<3$ it remains the only condition $\rho(\overset{\approx}{S})<4$).

\medskip

\noindent {\bf Conjecture 2.} {\it If $\overset{\approx}{S}$ is
tame, then $\rho(\overset{\approx}{S})\leq 4$ and
$\overline{\mu}(\overset{\approx}{S})\leq 4$.}

\medskip

We call a poset $\overset{\approx}{S}$ with a biequivalence
relation a {\it dyadic set} $D$ if $\dim \overset{\approx}{S}\leq
2$ and $\overset{\sim}{S}$ is chain ($s$ and $s^*$ are
comparable). This definition is equivalent to the definition given
in~\cite{4}. Let $\rho(D)=\rho(\overset{\approx}{S})$,
$\overset{\sim}{D}=\overset{\sim}{S}$.

We call the aggregate $K(D)$ {\it dyadic}. Not every aggregate of
dimension~2 is dyadic, but it is true (see \cite{15}),
if the aggregate is finitely represented.

Let $D_e=\{(x,y)\; |\; x,y\in D, x\Rightarrow y\}$. We call the elements
of $D_e$ {\it edges}. From i, ii it follows that
$(x,x^*)\not\in D_e$. If $\alpha=(x,y)\in D_e$, then
$(x,y)\approx(x^*,y^*)$ (i), and if $x<x^*$, then $y<y^*$ (if
$y^*<y$, then $x\lhd x^*\leq y^*\lhd y$ and $x\lhd y$). Set
$\alpha^*=(x^*,y^*)$.

The set $D_e$ is partially ordered. $(x,y)\,\overset{e}{\leq}\,
(\overline{x},\overline{y})$ if $\overline{x}\leq x$ and
$\overline{y}\leq y$. An edge is {\it maximal} (resp. {\it
short}), if one is maximal (resp. minimal) with respect to the order
relation. An edge is {\it long}, if it is not short
(for this, it does not have to be maximal).
If $x<y$, then we set $Eq(x,y)=D^\nc
\{x,y\}$, $eq(x,y)=\overset{\approx}{p} (Eq(x,y))$.
$\alpha\in D_e$ is {\it not equipped}, if $Eq(\alpha)=\varnothing$ and
{\it linearly equipped}, if $eq(\alpha)<\infty$. $D$ is {\it linearly equipped},
if every $\alpha\in D_e$ is linearly equipped.

$D'$ is a {\it $*$-subset} of a dyadic set $D$ if $a\in D'$
implies $a^*\in D'$. In this case $(D',\leq,\approx)$ is a dyadic
set and $\Rep D'\subset \Rep D$. A representation $T\subset \Rep
D$ is {\it faithful} if $T$ is not contained in $\Rep D'$ for any
$*$-subset $D'$.

$D$ is {\it faithful} if it has a faithful indecomposable
representation. If a dyadic set $\overline{D}$ is obtained from $D$
either by strengthening of the relation $\leq$, or by weakening of
the relation $\approx$ (a relation $\overline{R}$ is {\it weaker}
than a relation $R$ if $\overline{R}\subset R$), then
$\ind\overline{D}\subset \ind D$ (in particular,
$\ind\overset{\sim}{D}\subset \ind D$). $D$ is {\it critical} if
$|\ind D|=\infty$ and for any $\overline{D}$ and $D'$ as above
$|\ind \overline{D}|<\infty$ and $\ind D'|<\infty$.

An explicit list of critical dyadic sets would be long and
cumbersome (see~\cite{4}). Looking ahead note that
the results of~\cite{4}, which will be reformulated below,
imply that on a critical set a biequivalence is transitive
($a\Rightarrow b\Rightarrow c$ implies $a\Rightarrow c$, section~5);
we have no proof of this fact a priori.

In order to reformulate the criterion of finite representativity
of~\cite{4} we may assume that $|\ind\overset{\sim}{D}|<\infty$,
since $\ind\overset{\sim}{D}\subset \ind D$.
$\overset{\sim}{D}$ is a poset with equivalence for which the
criterion of finite representativity is given in proposition~5
($\rho(\overset{\sim}{D})<4$).

\begin{lemma}
If $|\ind\overset{\sim}{D}|<\infty$ and $a,a^*,b,b^*\in D$, then
either $a$ and $b$ are comparable, or $a^*$ and $b^*$ are comparable.
\end{lemma}

If $a\,\nc\, b$ and $a^*\,\nc\, b^*$, then it is not difficult to
see that none of the points $a$, $b$, $a^*$, $b^*$ is normal
$\overset{\sim}{p}(a)=\overset{\sim}{p}(b)=
\overset{\sim}{p}(a^*)=\overset{\sim}{p}(b^*)=\infty$,
$\rho(\overset{\sim}{D})\geq 4$, that contradicts to
proposition~5.

It is easy to see that $\rho(D)\geq \rho(\overset{\sim}{D})$. The
lemma below follows from the definitions and \cite{4}.

\begin{lemma} If $|\ind\overset{\approx}{D}|<\infty$, then
$\overset{\approx}{D}$ is linearly equipped. If $D$ is
linearly equipped and $\rho(\overset{\sim}{D})<4$, then $\rho(D)<4$.
\end{lemma}

A set $\{d_1\Rightarrow d_2\Rightarrow \cdots\Rightarrow
d_t\} (t\geq 1)\subset D$ is called {\it a strip} if any edge
$(d_i,d_{i+1})$ is short, there is no edge with the head $d_1$
(i.e. there is no $x\Rightarrow d_1$), and
there is no edge with the tail $d_t$
(i.e. there is no $d_t\Rightarrow y$).
For $t=1$ we have {\it an isolated} point $d$.

\begin{lemma}
If $\rho(\overset{\sim}{D})<4$, then there is no point, which is
the head or the tail of two short edges.
\end{lemma}

Let $a\Rightarrow b$ and $a\Rightarrow c$ are short edges.
$(a,b)\approx (a^*,b^*)$, $(a,c)\approx (a^*,c^*)$. If $b<c$, then
for $b\lhd c$, we have $a\lhd c$; and for $b\Rightarrow c$ the edge
$(a,c)$ is long. Analogously, cases $b>c$, $b^*<c^*$, $b^*>c^*$ are
excluded. So, $b\,\nc \, c$, $b^*\,\nc\, c^*$, that contradicts
to lemma~10.

\medskip

\noindent {\bf Corollary.} {\it If $\rho(\overset{\sim}{D})<4$,
then $\overset{\bullet}{D}$ is a disjoint union of strips.}

\medskip

A set $D$ is {\it bicomponent} if it can be decomposed
into ordinal sum of $D_1$ and $D_2$
($d_1\lhd d_2$ for $d_1\in D_1$, $d_2\in D_2$),
and each of $D_1$ and $D_2$ contains precisely one strip.

\begin{example} \rm $D_1=\begin{array}{c}
 \mbox{\includegraphics{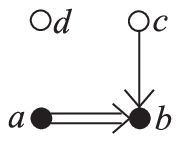}}
\end{array}$,  $D_2=\begin{array}{c}
 \mbox{\includegraphics{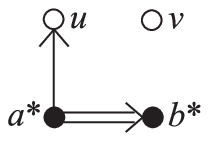}}
 \end{array}$.
 This set is finitely represented and has precisely one indecomposable representation.
 \end{example}

Representations of bicomponent sets are considered in~\cite{32}.
In a sense, representations of arbitrary finitely represented dyadic
sets can be reduced to representations of bicomponent sets.

A bicomponent set is {\it normal} if all big points of one of the
components $D_1$ and $D_2$ are normal. In example~4 bicomponent
$D$ is not normal. In~\cite{32} it is proved (and this is the main
result of the work) that a finitely represented faithful
bicomponent set is either normal, or one has the from given in
example~4.

For the formulation of a criterion of finite representativity of
dyadic sets, we need some more definitions.

Recall that $\mu(n_1,n_2,n_3)=n_1n_2+n_1n_3+n_2n_3+n_1n_2n_3$
(where $0\leq n_3\leq n_2\leq n_1$) (see lemma~6). In this section
we assume $\mu(\infty,0,0)=4$ and $\mu(\infty,n_2,n_3)=\infty$ for
$n_2\not=0$.

Let $\alpha=(x,y)\in D_e$. Set $\langle\alpha\rangle =\{z\in
D\;|\; x<z<y\}$ (then $(x,z),(z,y)\in D_e$)
$l(x,y)=l(\alpha)=|\langle \alpha\rangle|$. If $\rho(D)<4$, then
by lemma~12 $\omega(\langle \alpha\rangle)=1$.

If $u\in Eq(x,y)$, then $\{u\}\,\nc\, \langle \alpha\rangle$. Set
$Eq^-(\alpha)=\{v\in D\;|\; \{v\}\,\nc\,(\{y\}\cup\langle
\alpha\rangle)\}$, $Eq^+(\alpha)=\{v\in D\;|\; \{v\}\,\nc\,
(\{x\}\cup \langle \alpha\rangle)\}$.

A set $X\subset D_{\overset{\approx}{p}}$ is {\it bordering an
edge} $\alpha=(x,y)$ if $X=Z^-\cup Z^+\cup Z^e$, where $Z^-$,
$Z^+$, $Z^e$ are pairwise incomparable 1-chains, $Z^-\subset
Eq^-(\alpha)$, $Z^+\subset Eq^+(\alpha)$, $Z^e\subset Eq(\alpha)$,
and, moreover, if $l(\alpha)\geq 2$, then $Z^-=Z^+=\varnothing$
(a bordering set $X$, in general, is not uniquely defined by the edge $\alpha$).

If $X$ borders an edge $\sigma$, then set
$eq(\sigma,X)=\overset{\approx}{p}(Z^e)$,
$eq^*(\sigma,X)=eq(\sigma^*)+\min\{\overset{\approx}{p}(Z^-),
|2-l(\sigma)|\}+\min\{\overset{\approx}{p}(Z^+),|2-l(\sigma)|\}$\footnote{The
sign of absolute value at $2-l(\sigma)$ is missed in this formula in~\cite{4}.},
$\mu(\sigma,X)=\mu(eq(\sigma,X),eq^*(\sigma,X), l(\sigma))$.

\begin{example}\rm
$D=
\begin{array}{c}\mbox{\includegraphics{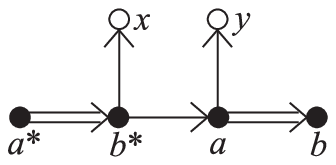}}\end{array}$,
$\sigma=(a,b)$, $X=\{x,y\}$, $\{x\}=Z^e$, $\{y\}=Z^-$,
$Z^+=\varnothing$, $eq(\sigma,X)=1$, $eq(\sigma^*)=0$,
$eq^*(\sigma,X)=1$, $l(\sigma)=0$, $\mu(\sigma,X)=1$.
\end{example}

The main theorem of \cite{4} (taking into account remark~7 and
lemma~11) can be reformulated in the following form.

\begin{proposition}
If $|\ind \overset{\sim}{D}|<\infty$, i.e.
$\rho(\overset{\sim}{D})<4$ (see proposition~5), then $|\ind
D|<\infty$ iff the following conditions are satisfied:

A. If $X$ borders $\sigma \in D_l$, then $\mu(\sigma,X)<4$.

B. If $\sigma=(a,b)\in D_e$, $l(\sigma)=1$, $eq(\sigma)=3$, then
$D^\nc (\{a\}\cup Eq(\sigma))= D^\nc (\{b\}\cup
Eq(\sigma))=\varnothing$.

C. If $\begin{array}{c}
 \mbox{\includegraphics{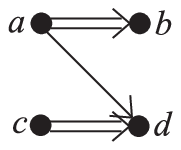}}
\end{array}\!\!\subset D$, then $\mu(eq(a,d),eq(a^*,b^*),eq(c^*,d^*))=0$.

\end{proposition}

(Linear equipment follows from A because if
$eq(\sigma^*)=\infty$, then $\mu(\sigma, X)\geq 4$.)

\begin{example} $D\supset Y=\begin{array}{c}
 \mbox{\includegraphics{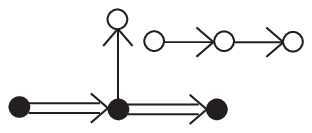}}
\end{array}$
\end{example}

This is a (unique, in some sense) example of non-fulfillment of
the condition~B. Consider this case in some more details. We have
$eq(\sigma)=3$ in the condition~B. At the same time, however, $|Eq(\sigma)|=3$
is possible (as in example~6, where $Eq(\sigma)$ consists of small points),
or $|Eq(\sigma)|=2$ ($Eq(\sigma)=\{s,t\}$, $s\in \overset{\bullet}{D}$,
$t\in \overset{\circ}{D}$, $D^\nc(s^*)=\varnothing$), or $Eq(\sigma)=y$
$(D^\nc(y^*)=z\in\overset{\circ}{D})$. Thus, the condition B corresponds to
several critical sets, however, they can be reduced to each other
since lemma~\cite{1} (see section~6) can be obviously generated
to posets with biequivalence if $x$ is an isolated (1-normal) point.

Call a critical set {\it primary} if it does not contain isolated
1-normal points (the rest can be reduced to them by this lemma).
Now we can say precisely that (up to duality) each primary set not
satisfying the condition~B contains $Y$ from example~6.

On the contrary, the condition~A excludes many (primary) critical sets
and can be a good example of an application of the function~$\mu$
(of course, here we may use $\rho$ instead of $\mu$).
First, suppose that $\overset{\approx}{S}$ is linearly equipped and
write out all cases when $\mu=4$ in the condition~A.
To shorten the number of cases, we assume
$eq^-=\min\{\overset{\approx}{p}(Z^-), |2-l(\sigma)|\}\geq eq^+=
\min\{\overset{\approx}{p}(Z^+),|2-l(\sigma)|\}$
(the rest cases are analogous).

 1) $l=0$,
$eq(\sigma)=eq(\sigma)=1$, $eq(\sigma^*)=4$ (or $eq(\sigma)=4$,
$eq(\sigma^*)=1$);

2) $l=0$, $eq(\sigma)=eq(\sigma)=2$, $eq(\sigma^*)=2$;

3) $l=0$, $eq(\sigma)=1$, $eq(\sigma^*)=3$, $eq^-=1$;

4) $l=0$, $eq(\sigma)=1$, $eq(\sigma^*)=2$, $eq^-=2$;

5) $l=0$, $eq(\sigma)=1$, $eq(\sigma^*)=2$, $eq^-=eq^+=1$;

6) $l=0$, $eq(\sigma)=eq(\sigma^*)=1$, $eq^-=2$, $eq^+=1$;

7) $l=0$, $eq(\sigma)=2$, $eq(\sigma^*)=1$, $eq^-=1$;

8) $l=0$, $eq(\sigma)=2$, $eq(\sigma^*)=0$, $eq^-=2$;

9) $l=0$, $eq(\sigma)=2$, $eq(\sigma^*)=0$, $eq^-=eq^+=1$;

10) $l=0$, $eq(\sigma)=4$, $eq^-=1$;

11) $l=1$, $eq(\sigma)=0$, $eq(\sigma^*)=4$ (or $l=1$,
$eq(\sigma)=4$, $eq(\sigma^*)=0$);

12) $l=1$, $eq(\sigma)=0$, $eq(\sigma^*)=3$, $eq^-=1$;

13) $l=1$, $eq(\sigma)=0$, $eq(\sigma^*)=2$, $eq^-=1$, $eq^+=1$;

14) $l=1$, $eq(\sigma)=1$, $eq(\sigma^*)=1$;

15) $l=1$, $eq(\sigma)=1$, $eq^-=1$;

16) $l=2$, $eq(\sigma)=2$, $eq(\sigma^*)=0$ (or $l=2$,
$eq(\sigma)=0$, $eq(\sigma^*)=0$);

17) $l=4$, $eq(\sigma)=1$, $eq(\sigma^*)=0$ (or $l=4$,
$eq(\sigma)=0$, $eq(\sigma^*)=1$).

In each of the listed cases we obtain a primary critical set
(see~\cite{4}) except the cases 8 and 15, where
corresponding $\overset{\approx}{S}$ are not critical
(in both cases, the set remains to be infinitely represented
after exclusion of the pair of equivalent points being the beginnings of strips).
From each primary set we can obtain several (sometimes one, but
always a finite number) critical sets similarly to the condition~B
and example~6.

Infinitely many critical sets correspond to the case
of a not linearly equipped $\overset{\approx}{S}$
$(eq(\sigma^*)=\infty)$.
Also infinitely many critical sets are obtained in case of
$\overset{\sim}{S}$ (see proposition~5).

Critical linearly equipped sets obtained from A, B contain one pair
of dual strips, whereas ones obtained from C contain two such pairs.

In other terms, the criterion of finite representativity of dyadic
sets is announced in~\cite{15} and proved in~\cite{26,4}.

Namely, for each dyadic set $D$ it is constructed a poset $C(D)$
so that $|\ind C(D)|<\infty$ is equivalent to $|\ind D|<\infty$.
Moreover, ``almost all'' representations of  $\ind D$ are {\it
multielementary} (see~\cite{31}), i.e. correspond to
representations of $\ind C(D)$. The unique faithful indecomposable
representation in example 4 is an example of a non-multielementary
representation. In some sense, the rest of non-multielementary
representations of finitely represented dyadic sets can be
obtained from this one, they are described in~\cite{27,28}.

Each of two formulations of the criterion of finite
representativity of dyadic sets (in terms of $C(D)$ and
proposition~7) evidently has its own advantages and shortages.
Proposition~7 is simpler for applications, closer to the explicit list
of critical sets (which is enough cumbersome) and allows to analyze
the qualitative structure of finitely represented dyadic sets
(the number of strips, transitivity and so on).

From \cite{24} it follows that for $\dim{\mathcal A}=3$ and
$|\ind {\mathcal A}|<\infty$ a representation of the aggregate ${\mathcal A}$
is a representation of {\it a triadic} set
\cite{24,5}. A criterion of their finite
representativity generalizing proposition~7 is formulated and
proved to one side in~\cite{5}. There the function $\mu$ (or $\rho$)
again plays an important role;
$\mu$ appears there also in one case, which is essentially
different from the considered ones in this article.

\section{Criteria of finite representativity and tameness\\
for semilinear marked quivers}

In this section we give criteria of finite representativity and
tameness for some types of marked quivers in terms of functions
$\rho$ on the base of the results of \cite{9} and section~3. Also
we show that the problem of finite representativity and tameness
for arbitrary $k$-marked quivers can be reduced to analogous
problems for representations of aggregates.

An aggregate $K(S)$ is {\it linear} if  $w(S)=1$.
An aggregate $K(\overset{\sim}{S})$ and a poset with equivalence
$\overset{\sim}{S}$ are {\it semilinear} if
$\dim\overset{\sim}{S}\leq 2$, each element $x\in S$ is comparable
with no more than one element $y$, and if such $y$ exists,
then $\dim x=1$ (here we consider a linear aggregate as a partial case
of semilinear one). A quiver $\overline Q_v$ is {\it semilinearly marked}
if for each $x\in Q_v$ the subaggregate $K(x)$ is semilinear. Denote
the aggregate $K(S,\leq)$ by $K^n$ if $w(S)=1$, $|S|=n$;
we assume that $K^0$ is the ``empty'' aggregate ($\Ob K^0=\varnothing$).

If $\A$ and $\B$ are two aggregate, then by $\A\oplus \B$ we denote
the aggregate defined via $\Ob(\A\oplus\B)=\{A\oplus B|A\in\Ob \A, B\in
\Ob \B\}$, $\Hom (A_1\oplus B_1,A_2\oplus B_2)= \A(A_1,A_2)\oplus
\B(B_1,B_2)$.

For an aggregate $\A$, its dual aggregate $\A^\circ$
can be naturally defined: $\Ob
\A^\circ=\{A^*=\Hom_k(A,k)|A\in \Ob \A\}$,
$\A^\circ(A^*,B^*)=\{\varphi^*| \varphi\in \A(B,A)\}$\footnote{We
have no proof that aggregates $\A$ and $\A^\circ$ are simultaneously finitely
represented or tame, though it seems completely obvious.}.

To a semilinear marked quiver $\overline Q$, we attach a $v$-graph
$\overline{\Gamma}=\overline{\Gamma}(Q)$ (see section~3) assuming
$v(x)=\dim K(x)$ if $K(x)$ is linear and $v(x)=\infty$ otherwise.

The main theorem of~\cite{9} and proposition~3 (taking into account remark~4)
imply the following statement.

\begin{proposition} Semilinear marked quivers $\overline Q$
is finitely represented (resp. tame) if and only if $\rho$-degree of each
vertex of $\overline{\Gamma}(Q)$ is less than 4, (resp. $\max
g_\rho(x)=4$, where $x\in Q_v$).
\end{proposition}

If $\overline Q$ is an arbitrary marked quiver, then preserve the
sense of $v(x)$ if $K(x)$ is semilinearly marked and write $v(x)>\infty$
in the opposite case.

\begin{lemma}[\cite{9}]
If $v(x)>\infty$ ($x\in Q_v$), then $Q$ is finitely represented
(resp. tame) if and only if
\begin{enumerate}
\item
The graph $\Gamma(Q)$ is $A_l$ ($l\geq 2$), $x=a_1$.
\item
$v(a_i)=1$ for $1<i<l$, $v(a_l)<\infty$.
\item
The aggregate $\A\oplus D$ is finitely represented (resp. tame), where
$\A=K(x)$ or $K(x)^\circ$ if $x$ is the head or the tail of an arrow
of quiver~$Q$ respectively, $\D=K^{t}$, $t=l+v(a_l)-3$.
\end{enumerate}
\end{lemma}

Criteria of finite representativity and tameness for quivers
non-semilinearly marked by $K(\overset{\sim}{S})$ or $K(D)$ can
be reduced to representations of $\overset{\sim}{S}$ and $D$ itself
by lemma~13. Write out these criteria in evident forms.

Let $x\in Q$, $v(x)>\infty$, let conditions 1 and 2 of lemma~13
be satisfied, $t=|Q_v|+v(a_l)-3$, and let $K(x)$ be
$K(\overset{\sim}{S})$ or $K(D)$. First, suppose that
$\overline{Q}$ is a marked quiver, where
$K(x)=K(\overset{\sim}{S})$. Then $\overline{Q}$ is finitely
represented in the following cases:
\begin{enumerate}
\item[1)] $t=1$; $\dim \overset{\sim}{S}\leq 2$, $\rho(\overset{\sim}{S})<3$;
\item[2)] $t\in\{2,3,4\}$; $\dim \overset{\sim}{S}\leq 2$,
$\rho(\overset{\sim}{S})<3-\frac{t-1}{t+1}$ (if $t=4$, then
$\rho(\overset{\sim}{S})<2,4$; $S\not\supset \widehat N$, i.e $S$
is an ordinal sum of posets of the form (1), (1,1) and (1,2) (see
remark 3), moreover, if $\dim s=2$, then $s\not\in\{1,2)$, and if
$s\in(1,1)$, then $S^\nc(s^*)=\varnothing$).
\end{enumerate}

$\overline{Q}$ is tame in the following cases:
\begin{enumerate}
\item[1)${}'$] $t=1$;
$\dim \overset{\sim}{S} \leq 3$, and if $\dim \overset{\sim}{S} <
3$, then $\rho(\overset{\sim}{S})=3$, and if $\dim s= 3$ $(s\in
S)$, then $S^\nc(s)=\varnothing$ and $\rho(S)\leq 3$;

\item[2)${}'$] $t\in\{2,3,4,5\}$; $\dim \overset{\sim}{S}
\leq 2$, $\rho(\overset{\sim}{S})=3 -\frac{t-1}{t+1}$ (if $t=5$,
then $\rho(\overset{\sim}{S})=2\frac 13$ and $\overset{\sim}{S}$
has the same form as in 2) for $t=4$).
\end{enumerate}

Let $Q(x)=K(D)$ and $D\not= \overset{\sim}{D}$. The corresponding
marked quiver is finitely represented iff  $t=1$, $\rho(D)<3$, all
edges of $D$ are short and not equipped and if for some
$(\sigma,X)$, $Z^-\not=\varnothing$ and $Z^+\not=\varnothing$,
then $\overset{\approx}{p}(Z^-)=\overset{\approx}{p}(Z^+)=1$.
A quiver $Q$ is marked by triadic sets can be finitely represented
only in case  $Q=\Delta$ and if one of the points is trivially marked,
i.e. $\Rep Q_M$ is equivalent to the category of
representations of a triadic set (marking the second vertex or its dual one).

The authors are very grateful to I.R. Shafarefich and to participants
of his seminar for the remarks essentially used for writing this article.

\end{document}